%% file: main.tex
\documentclass[journal]{IEEEtran}
\ifCLASSINFOpdf
   \usepackage[pdftex]{graphicx}
   \graphicspath{{figures/}}
\else
\fi
\usepackage{amsmath}
\usepackage{amsfonts} 
\newcommand{\specialcell}[2][c]{%
  \begin{tabular}[#1]{@{}c@{}}#2\end{tabular}}

\usepackage[export]{adjustbox}
\usepackage{enumitem}
\usepackage{comment}
\usepackage{booktabs}
\usepackage{color}
\usepackage[table]{xcolor}
\usepackage{arydshln}
\usepackage{float}
\usepackage{multirow}

\usepackage[linesnumbered,ruled]{algorithm2e}
\SetAlFnt{\small}
\SetAlCapFnt{\small}
\SetAlCapNameFnt{\small}

\SetCommentSty{mycommfont}

\begin{document}
%
% paper title
% Titles are generally capitalized except for words such as a, an, and, as,
% at, but, by, for, in, nor, of, on, or, the, to and up, which are usually
% not capitalized unless they are the first or last word of the title.
% Linebreaks \\ can be used within to get better formatting as desired.
% Do not put math or special symbols in the title.
\title{GPSAF: A Generalized Probabilistic Surrogate-Assisted Framework for Constrained Single- and Multi-objective Optimization}

% author names and IEEE memberships
% note positions of commas and nonbreaking spaces ( ~ ) LaTeX will not break
% a structure at a ~ so this keeps an author's name from being broken across
% two lines.
% use \thanks{} to gain access to the first footnote area
% a separate \thanks must be used for each paragraph as LaTeX2e's \thanks
% was not built to handle multiple paragraphs
%

\author{Julian~Blank,~\IEEEmembership{Student Member,~IEEE} and
        Kalyanmoy~Deb,~\IEEEmembership{Fellow,~IEEE}
\thanks{Julian Blank and Kalyanmoy Deb are at Michigan State University, East Lansing, USA (email: \{blankjul,kdeb\}@msu.edu).}
}

% note the % following the last \IEEEmembership and also \thanks -
% these prevent an unwanted space from occurring between the last author name
% and the end of the author line. i.e., if you had this:
%
% \author{....lastname \thanks{...} \thanks{...} }
%                     ^------------^------------^----Do not want these spaces!
%
% a space would be appended to the last name and could cause every name on that
% line to be shifted left slightly. This is one of those "LaTeX things". For
% instance, "\textbf{A} \textbf{B}" will typeset as "A B" not "AB". To get
% "AB" then you have to do: "\textbf{A}\textbf{B}"
% \thanks is no different in this regard, so shield the last } of each \thanks
% that ends a line with a % and do not let a space in before the next \thanks.
% Spaces after \IEEEmembership other than the last one are OK (and needed) as
% you are supposed to have spaces between the names. For what it is worth,
% this is a minor point as most people would not even notice if the said evil
% space somehow managed to creep in.

% The paper headers
\markboth{}%
{Shell \MakeLowercase{\textit{et al.}}: Bare Demo of IEEEtran.cls for IEEE Journals}
% The only time the second header will appear is for the odd numbered pages
% after the title page when using the twoside option.
%
% *** Note that you probably will NOT want to include the author's ***
% *** name in the headers of peer review papers.                   ***
% You can use \ifCLASSOPTIONpeerreview for conditional compilation here if
% you desire.

% If you want to put a publisher's ID mark on the page you can do it like
% this:
%\IEEEpubid{0000--0000/00\$00.00~\copyright~2015 IEEE}
% Remember, if you use this you must call \IEEEpubidadjcol in the second
% column for its text to clear the IEEEpubid mark.

% use for special paper notices
%\IEEEspecialpapernotice{(Invited Paper)}

% make the title area
\maketitle

% As a general rule, do not put math, special symbols, or citations
% in the abstract or keywords.
\begin{abstract}
Significant effort has been made to solve computationally expensive optimization problems in the past two decades, and various optimization methods incorporating surrogates into optimization have been proposed.
Most research focuses on either exploiting the surrogate by defining a utility optimization problem or customizing an existing optimization method to use one or multiple approximation models.
However, only a little attention has been paid to generic concepts applicable to different types of algorithms and optimization problems simultaneously. 
Thus this paper proposes a generalized probabilistic surrogate-assisted framework (GPSAF), applicable to a broad category of unconstrained and constrained, single- and multi-objective optimization algorithms.
The idea is based on a surrogate assisting an existing optimization method.
The assistance is based on two distinct phases, one facilitating exploration and another exploiting the surrogates. 
The exploration and exploitation of surrogates are automatically balanced by performing a probabilistic knockout tournament among different clusters of solutions.
A study of multiple well-known population-based optimization algorithms is conducted with and without the proposed surrogate assistance on single- and multi-objective optimization problems with a maximum solution evaluation budget of 300 or less. The results indicate the effectiveness of applying GPSAF to an optimization algorithm and the competitiveness with other surrogate-assisted algorithms. 
\end{abstract}

% Note that keywords are not normally used for peer review papers.

% For peer review papers, you can put extra information on the cover
% page as needed:
% \ifCLASSOPTIONpeerreview
% \begin{center} \bfseries EDICS Category: 3-BBND \end{center}
% \fi
%
% For peer review papers, this IEEEtran command inserts a page break and
% creates the second title. It will be ignored for other modes.
\IEEEpeerreviewmaketitle

\section{Introduction}

% expensive problems to surrogate methods
\IEEEPARstart{M}{any} optimization problems are computationally expensive and require the execution of one or multiple time-consuming objective and constraint functions to evaluate a solution. Expensive optimization problems (EOPs) are especially important in practice and are omnipresent in all kinds of research and application areas, for instance,~Agriculture~\cite{2019-roy-crop-yield},~Engineering~\cite{2016-yin-crashwrothiness-design}, Health Care~\cite{2016-lucidi-application-health}, or Computer Science~\cite{2019-lu-nsga-net}. Often the expensive solution evaluations (ESEs) are caused by running a simulation, such as Computational Fluid Dynamic~\cite{1995-anderson-cfd}, Finite Element Analysis~\cite{1991-szabo-fea}, or processing a large amount of data~\cite{2019-luo-data-driven,2020-wang-data-driven-forest}.
Most of such these simulation-based or data-intensive ESEs are black-box in nature~\cite{2002-olafsson-simopt} and no gradient information is available or even more time-consuming to derive.
Thus, the optimization method must be designed for a significantly limited evaluation budget without any assumptions about the problem's fitness landscape.

The majority of methods proposed for solving EOPs incorporate so-called surrogate models (or metamodels)~\cite{2011-jin-surrogate-recent} into the optimization procedure. The surrogate model provides an approximate solution evaluation (ASE) with less computational expense to improve the convergence behavior.
A well-known research direction with the significant effort being made in the past is referred to as "Efficient Global Optimization" (EGO)~\cite{1998-jones-ego}.
In EGO, solutions are evaluated in each iteration based on the optimum of a utility optimization problem -- also known as infill criterion -- commonly defined by the surrogate's value and error predictions from Kriging~\cite{1951-krige-mine}.
The method can be summarized as a \emph{fit-define-optimize} procedure: A surrogate model is fitted, a utility problem based on the surrogate predictions is defined, and then optimized to obtain an infill solution. 
Original limitations such as the evaluation of a single point per iteration, the lack of constraint handling, or dealing with multiple objectives have been investigated, and extensions have been proposed~\cite{2016-haftka-paralllel,2006-knowles-parego}.
Nevertheless, adding such additional requirements further complicates the utility problem and makes it significantly more challenging.
The surrogate's role in such \emph{fit-define-optimize}  methods is \emph{critical} because of the direct dependency on the infill criteria defined based on approximations and error predictions.

Another research direction pursued by researchers is the organic incorporation of surrogates into an existing optimization ~\cite{2002-jin-approximate-fitness}. Such approaches aim to improve the convergence behavior of a ``baseline'' algorithm and, thus, the anytime performance.
Researchers have explored different ways of incorporating surrogates into well-known population-based algorithms, such as genetic algorithms (GA)~\cite{1989-goldberg-ga}, differential evolution (DE)~\cite{1996-storn-de}, or particle swarm optimization (PSO)~\cite{1995-kennedy-pso}.
All surrogate-assisted algorithms must find a reasonable trade-off between exploiting the knowledge of the surrogate and still exploring the search space.
On the one hand, researchers have investigated methods adding a surrogate with lighter influence on the original algorithm, for instance, using surrogate-based pre-selection in evolutionary strategy~\cite{2002-emmerich-preselection} or a predictor for the individual replacement in DE~\cite{2011-lu-de-classification}.
On the other hand, the behavior of an existing algorithm might almost entirely rely on the surrogate predictions and guide the search significantly.
For instance, a global and local surrogate have been incorporated into PSO to solve expensive large-scale optimization problems~\cite{2019-wang-local-global} or into DE for expensive constrained optimization problems~\cite{2017-sun-SACOSO}.
The existence of numerous variants of surrogate-assisted algorithms indicates that many different ways of using surrogates during optimization method exist, but also that no best practice procedure has been established yet~\cite{2020-stork-issues}.
\begin{figure}[t]
    \centering
    % trim=left bottom right top, clip. (add fbox)
    \includegraphics[page=4,width=0.95\linewidth,trim=0.2cm 9.5cm 16.2cm 0.5cm,clip]{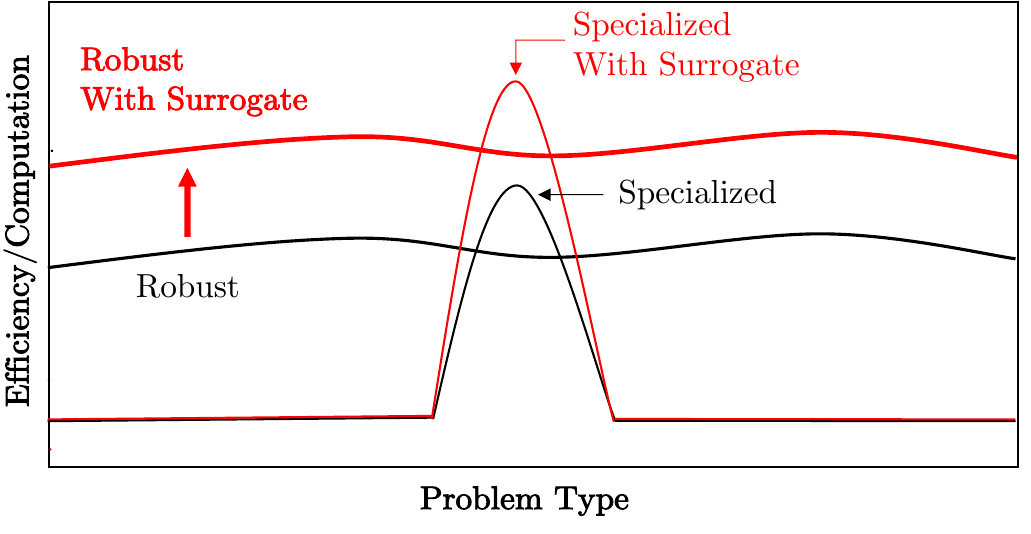}
    \caption{Robustly Adding Surrogate-Assistance to Population-Based Algorithms (Illustration inspired from~\cite{1989-goldberg-ga}).}
 \label{fig:surrogate}
 \vspace{-3mm}
 \end{figure}

The need for more generalizable concepts in surrogate-assisted optimization has been identified, and frameworks aiming to solve a broader category of optimization problems have been proposed~\cite{2010-lim-generalizing,2020-cai-generalied-high-dim}. Existing frameworks provide a generic method for different problems types using the \emph{fit-define-optimize} strategy or by replacing locally optimized solutions (on ASEs) before their computationally expensive evaluation. 
One shortcoming of existing frameworks for solving EOPs is their design being primarily algorithm-dependent and their limitations to transferring it to different optimization methods.
Thus, this matter shall be addressed in this paper by proposing a novel generalized probabilistic surrogate-assisted framework (GPSAF), adding surrogate assistance to population-based algorithms.
In contrast to other surrogate-assisted algorithms that customize a \emph{specific} algorithm, our goal is to provide a scheme to add surrogate assistance to a whole algorithm class. 
Figure~\ref{fig:surrogate} -- inspired from 
~\cite{1989-goldberg-ga} -- shows different problem types on the $x$ and the efficiency of algorithms on the $y$-axis. Whereas specific algorithms can be customized and thus be a specialist for a specific problem type, this study considers algorithms that can solve a broad category of problem types and adds surrogate assistance to them.
Even though specialized surrogate-assisted algorithms are likely to outperform a generic concept on a specific problem type, the merit of this study is its broad applicability.
The main contributions of this paper can be summarized as follows:

\begin{enumerate}[label=(\roman*)]
\item We provide a categorization of existing surrogate-assisted algorithms regarding their surrogate usage. Methods are distinguished based on their surrogate's impact and importance during optimization and identify what has been less attention paid to in the past.

\item We propose a framework that applies to all kinds of population-based optimization algorithms. The framework enables existing optimization methods designed for unconstrained or constrained, single or multi-objective optimization to become assisted by a surrogate. Intuitive hyper-parameters can control the surrogate's impact. A specific setting to even disable the surrogate usage entirely exists, which demonstrates the truly surrogate assistant behavior.

\item We propose a novel way of dealing with surrogate prediction error algorithmically. 
In contrast to existing surrogate-assisted methods, we are using the search pattern of the search of an algorithm on the surrogate instead of only using final solutions. 
The prediction error of the surrogates and the search space exploration is addressed using a probabilistic knockout tournament selection. The surrogate prediction error is incorporated into the tournament selection and reliably balances the exploitation-exploration trade-off based on the surrogate's accuracy.

\end{enumerate}

% Structure of this paper
In the remainder of this paper, we first discuss related work of surrogate-assisted optimization and its challenges in~Section~\ref{sec:background} before proposing the generalized probabilistic surrogate-assisted framework (GPSAF) in Section~\ref{sec:method}.
A comparison of GPSAF and other state-of-the-art surrogate-based methods is provided in~Section~\ref{sec:results}.
Finally, conclusions are drawn, and future work is presented in Section~\ref{sec:conclusion}.

\section{Background}
\label{sec:background}

\begin{figure*}
    \centering
    % trim=left bottom right top, clip. (add fbox)
    \includegraphics[page=1,width=0.85\linewidth,trim=0 9.5cm 0cm 0,clip]{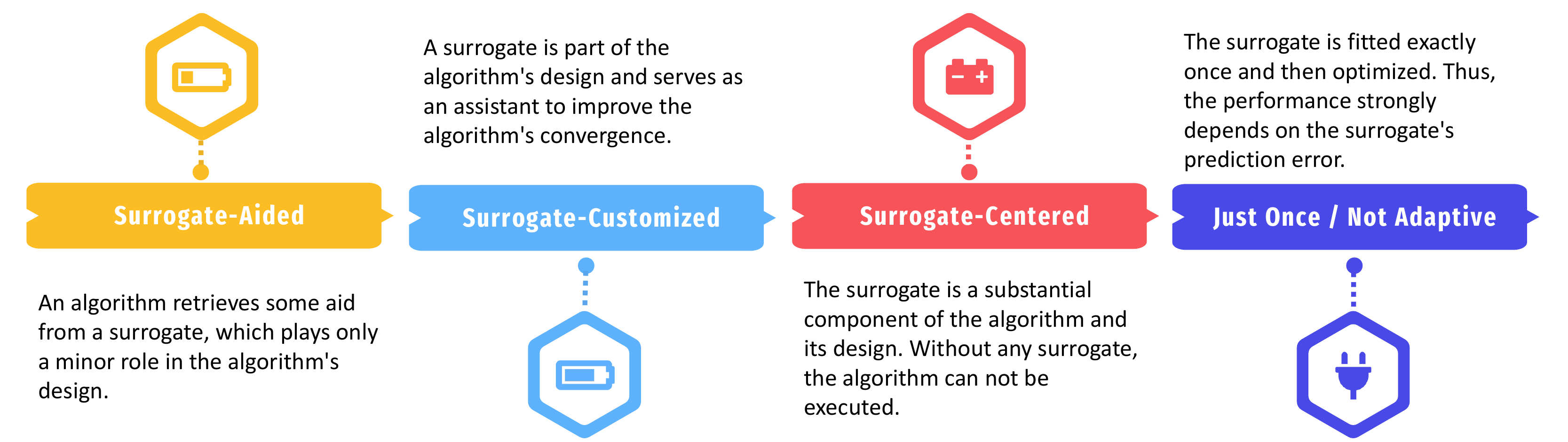}
    \caption{Different Roles of Surrogates in the Design of an Algorithm.}
 \label{fig:overview}
 \vspace{-3mm}
 \end{figure*}

In this section, the brief overview of methods in the previous section shall be enriched with details, and surrogate-assisted algorithms are categorized regarding their surrogate incorporation.

Surrogate-assisted methods can be roughly put into one of the following categories based on the surrogate's involvement: aided, customized, centered, or once (see Figure~\ref{fig:overview}). The latter describes the early development of optimization using an approximation model, fitted exactly once during optimization and never updated (once).
Algorithms that perform an update of the surrogate can mostly depend on its predictions (centered) or use it as an assistant in an existing method to improve the convergence behavior (aided, customized). The surrogate's role and dependency on the algorithm's design are vital for generalization and, thus, shall be spent special attention to.
Next, a thematic overview of these different types of surrogate involvements in an algorithm is given.

\begin{table}[b]
\renewcommand{\arraystretch}{1.2}
\label{tbl:category}
\caption{Categorization regarding the Surrogate's Role in an Optimization Algorithm.}
\centering
\begin{tabular}{p{1.3cm} | p{6cm}}
\toprule
\textbf{Category} & \textbf{Algorithm / Study} \\ \midrule
 \textbf{Aided} & MAES~\cite{2002-emmerich-preselection}, SVC-DE~\cite{2011-lu-de-classification}\\ \hline
 \textbf{Customized} &
 MOEAD-EGO~\cite{2010-zhang-moead-ego},
 K-RVEA~\cite{2016-chugh-krvea},
 HSMEA~\cite{2019-habib-hsmea},
 CSEA~\cite{2019-pan-csea},
 PAL-SAPSO~\cite{2019-lv-surr-pso},
 CAL-SAPSO~\cite{2017-wang-pso-committee-active-learning}\\ \hline
 \textbf{Centered} &  EGO~\cite{1998-jones-ego},
 ParEGO~\cite{2006-knowles-parego},
 SMS-EGO~\cite{2008-ponweiser-sms-ego},
 Max-Min~SAEA~\cite{2006-ong-max-min},
 SACOBRA~\cite{2017-bagheri-sacobra},
 SABLA~\cite{2017-islam-sabla-bilevel},
 GS-MOMA~\cite{2010-lim-generalizing},
 GSGA~\cite{2020-cai-generalied-high-dim},
 MISO~\cite{2016-mueller-miso},
 GOSAC~\cite{2017-mueller-gosac}
 \\ \hline
 \textbf{Just Once} &  \cite{2014-lei-level-design}, \cite{2013-duan-emd}, \cite{2006-dorica-electromagnitic-neural} \\ \hline
\end{tabular}

\end{table}

Especially in the early phase of surrogate-based optimization, the surrogate was fitted \emph{only once} and optimized. Thus, the optimization's outcome entirely depends on the accuracy of the surrogate model, assuming an efficient optimizer is used to obtain the surrogate's optimum.

% dependent
The limitation of fitting a surrogate only once has soon been overcome by a more adaptive approach known as EGO (Efficient Global Optimization)~\cite{1998-jones-ego}.
Kriging~\cite{1951-krige-mine} is used as a surrogate and provides predictions as well as a measure of uncertainty for each point. The prediction and uncertainty together define the so-called acquisition function (or infill criterion), such as the expected improvement~\cite{1998-jones-ego} or probability of improvement~\cite{2001-jones-tax-ego} aiming to balance exploitation and exploration simultaneously. The optimization of the acquisition function results in an infill solution, which is first evaluated and then added to the model. The procedure is repeated until a termination criterion is met. The limitation of finding only a single new solution in each iteration has been investigated thoroughly, and multi-point EGO approaches have been proposed~\cite{2010-viana-ego-parallel,2019-beaucaire-ego-multi-point,2020-berveglieri-parallel}.
Moreover, the concept has been generalized to solve multi-objective optimization problems by using decomposition~\cite{2006-knowles-parego,2010-zhang-moead-ego} or replacing the objective with a performance indicator based metric~\cite{2008-ponweiser-sms-ego}. The idea has also been extended to handle constraints, which is especially important for solving real-world optimization problems~\cite{2017-bagheri-constr-ego}.
Instead of using acquisition functions to address the surrogate's uncertainty, algorithms based on trust regions have been proposed.
Inevitable, updating the trust-region radii becomes vital for the algorithm's performance~\cite{2003-ong-evolution-expensive}. Whereas original studies were limited to unconstrained single-objective optimization, the surrogate-assisted trust-region concept has been generalized to constrained and also multi-objective optimization~\cite{2010-lim-generalizing,2020-cai-generalied-high-dim}.
Apart from the approaches discussed above, the direct usage of surrogates in an algorithm has been explored in various areas, for instance, bi-level optimization~\cite{2006-ong-max-min,2017-islam-sabla-bilevel} or mixed-integer optimization~\cite{2016-mueller-miso}.
All these approaches have in common that the algorithm has been designed with a strong dependency on the surrogate model. Thus, the surrogate's suitability and accuracy are critical for the optimization's success. Inaccurate surrogate predictions and error estimations, inevitably occurring in large-scale optimization problems, are known to be problematic~\cite{2020-stork-issues}.

% surrogate-assisted
In contrast to algorithms being designed based on surrogates, researchers have investigated surrogates' incorporation into existing optimization methods. Such approaches are also known as \emph{surrogate-assisted} algorithms, emphasizing the surrogate's role as an assistant during optimization.
In our categorization, surrogate-assisted algorithms are split up into two categories. On the one hand, algorithms can be aided by a surrogate where only minor changes of the original algorithm design are made; on the other hand, surrogate-customized methods where the algorithm has a significant impact on the algorithm's design.
Because the judgment of \emph{impact} is subjective, the transition between both categories is somewhat fluent.

The benefit of surrogate-aided algorithms is that with relatively minor modifications, a surrogate has been incorporated, and the performance has been improved~\cite{2011-jin-surrogate-recent}. One well-known approach is a pre-selection (or pre-filtering) which uses a surrogate to select a subset of solutions that usually would be evaluated on the EOP~\cite{2002-emmerich-preselection}.
Moreover, instead of changing the behavior in a generation, surrogates have also been used across generations by switching between the expensive evaluation and surrogate predictions entirely for some iterations~\cite{2002-jin-approximate-fitness,1998-ratle-surr-ga}.
Another example for a surrogate-influenced algorithm is modifying a memetic algorithm (genetic algorithm with local search) by executing the usually evaluation-intensive local search on the surrogate~\cite{2007-zhou-local-vs-global}.

Besides surrogate-assisted methods with relatively minor modifications of existing algorithms, optimization methods can be customized to incorporate surrogate usage.
This let arise surrogate-assisted variants of well-known algorithms, such as lqCMAES~\cite{2019-hansen-lqcmaes} derived from CMAES~\cite{2001-hansen-cmaes}, KRVEA~\cite{2016-chugh-krvea} and HSMEA~\cite{2019-habib-hsmea} based on RVEA~\cite{2016-cheng-rvea}, MOEAD-EGO~\cite{2010-zhang-moead-ego} as an improvement of MOEAD~\cite{2007-zhang-moead}, or CAL-PSO~\cite{2017-wang-pso-committee-active-learning} based on PSO~\cite{1995-kennedy-pso} to name a few. Each surrogate-assisted variant is in principle based on an algorithm originally developed for computationally more inexpensive optimization problems but customizes the default behavior, for instance, by one or multiple local or global surrogates, implementing a sophisticated selection after surrogate-based optimization.

The increasing number of surrogate-assisted algorithms shows its importance and relevance in practice. Indisputably, approaches and ideas directly designed for and centered on one or multiple surrogates have their legitimacy but are somewhat tricky to use for newly proposed algorithms.
The existence of many different surrogate-based algorithms also indicates the absence of a best practice procedure and the need for more generic methods. Thus, this study shall address this research gap and propose a surrogate-assisted framework of algorithms applicable to a broad category of optimization methods.

\section{Methodology}
\label{sec:method}

% interface and assumptions
One of the major challenges when proposing a generalized optimization framework is the number and strictness of assumptions being made. 
On the one hand, too many assumptions restrict the applicability; on the other hand, too few assumptions limit the usage of existing elements in algorithms.
In this study, we target any type of population-based algorithm with two phases in an iteration: the process of generating new solutions to be evaluated (infill) and a method processing evaluated infill solutions (advance). 
With these two methods, running an algorithm can be summarized by the pseudo-code shown in Algorithm~\ref{alg:infill-advance}. Until the algorithm $\Phi$ has been terminated, the \emph{infill} method returns a set of new designs $X$ to be evaluated. After obtaining the objective $F$ and constraint $G$ values for design, the algorithm is advanced by providing the evaluated solutions \{$X$, $F$, $G$\}. 
By looking at this interface, we further make two (weak) assumptions. First, we do not assume that $X$ needs to be identical with the suggested designs from \emph{infill} (Line 2 and 4), but can also be modified. Second, the \emph{infill} method is non-deterministic, resulting in different designs $X$ whenever called. Both assumptions can be considered weak because most population-based algorithms already fulfill them.
So, how can existing optimization methods be described into \emph{infill} and \emph{advance} phases?
Genetic algorithms (GAs) generate new solutions using evolutionary recombination-mutation operators and then process them using an environmental survival selection~\cite{1989-goldberg-ga} operator; PSO methods create new solutions based on a particles' current velocity, personal best, and global best, and process the solutions using a replacement strategy~\cite{1995-kennedy-pso}; CMAES samples new solutions from a normal distribution, which is then updated in each iteration~\cite{2001-hansen-cmaes}. Shown by well-known state-of-the-art algorithms following or being suitable to be implemented in this optimization method design pattern, this seems to be a reasonable assumption to be made for a generic framework. Moreover, it is worth noting that some researchers and practitioners also refer to the pattern as \emph{ask-and-tell} interface.

{\small
% \linespread{1.0}\selectfont
\begin{algorithm}[t]
\SetKwInput{Input}{Input~}
\Input{Algorithm $\Phi$
}

\vspace{1mm}

\While{$\Phi$ \texttt{has not terminated}}{

$X \gets \Phi.\tt{infill}()$

$F, G \gets \tt{evaluate}(X)$

$\Phi.\tt{advance}(X, F, G)$

}

\caption{Infill-And-Advance Interface}
\label{alg:infill-advance}
\end{algorithm}}

{\small
% \linespread{1.0}\selectfont
\begin{algorithm}[t]
\SetKwInput{Input}{Input~}
\Input{Algorithm $\Phi$, Surrogate Tournament Pressure $\alpha$  ($\geq 1$), Number of Simulated Iterations $\beta$ ($\geq 0$), Replacement Probability Exponent $\gamma$,
Maximum Number of Solution Evaluations $\tt{SE}^{(\max)}$
}

\vspace{1mm}

\tcc{Sample Design of Experiments (DOE)}

$\tt{A} \gets \emptyset ; \, \tt{P} \gets \emptyset ; \,            \tt{Q} \gets \emptyset; \,              \tt{U} \gets \emptyset ; \,             \tt{e} \gets \emptyset $

$\tt{A.X} \gets \tt{doe()}; \quad \tt{A.F, A.G} \gets \tt{evaluate(A.X)}$
\label{alg:gpsaf-doe}

\vspace{1mm}

\While{$\tt{size}(A) < \tt{SE}^{(\max)}$}{
\label{alg:gpsaf-while}

    \vspace{1mm}
    \tcc{Infill sols. from baseline algorithm}
    $\tt{P.X} \gets \Phi\tt{.infill()}$
    \label{alg:gpsaf-default-infill}

    \vspace{1mm}
    \tcc{Estimate error - only initially}
    \lIf{$\tt{e} = \emptyset$}{$\tt{e} \gets \tt{estm\_error(A.X, A.F, A.G)}$}

    \vspace{1mm}
    \tcc{Surrogates for each obj. and constr.}
    $S \gets \tt{fit}(\tt{A.X, A.F, A.G})$
    \label{alg:gpsaf-fit-surr}
    
    $\tt{P.}\hat{F}, \, \tt{P.}\hat{G} \gets S.\tt{predict}(\tt{P.X})$
    \label{alg:gpsaf-surr}

    \vspace{1mm}
    \tcc{Surrogate Influence $(\alpha)$}

    \ForEach{$k \gets 2$ \KwTo $\alpha$}{
    \label{alg:gpsaf-start_alpha}

        $\tt{Q.}X \gets \Phi\tt{.infill()}$

        $\tt{Q.}\hat{F}, \tt{Q.}\hat{G} \gets S.\tt{predict}(\tt{Q.X})$

        \ForEach{$j\gets1$ \KwTo  $\tt{size(Q)}$}{
            \lIf{\textbf{not} $\tt{dominates(P[j], Q[j])}$}{
            \label{gpsaf:alpha-comp}
              $\tt{P[j] = Q[j]}$
            }
        }

    }
    \label{alg:gpsaf-end_alpha}

    \vspace{1mm}
    \tcc{Surrogate Bias $(\beta)$}
    $\Phi^{\prime} \gets  \tt{copy}(\Phi)$
    \label{alg:gpsaf-copy}
    \label{alg:gpsaf-start_beta}

    $U \gets \emptyset$
    
    \ForEach{$k \gets 1$ \KwTo $\beta$}{
    
        $\tt{Q.}k \gets k$ 
    
        $\tt{Q.}X \gets \Phi^{\prime}.\tt{infill()}$
        \label{alg:gpsaf-beta_infill}
        
        $\tt{Q.}\hat{F}, \, \tt{Q.}\hat{G} \gets S.\tt{predict(Q.X)}$
        \label{alg:gpsaf-beta_predict}

        \ForEach{$j\gets1$ \KwTo  \tt{size(Q)}}{
            $i\gets \tt{closest}(\tt{P.}X, \tt{Q[j].}X)$
            
            $\tt{U[i]} \gets \tt{U[i]} \cup \tt{Q[j]}$
        }

        $\Phi^{\prime}.\tt{advance(Q.X, Q.}\hat{F}\tt{, Q.}\hat{G})$
        \label{alg:gpsaf-beta_advance}

    }
    
    $\tt{V} \gets \tt{list()}$
    
    \ForEach{$j\gets1$ \KwTo $\tt{size(U)}$}{
        $\tt{V} \gets \tt{V} \cup \tt{prob\_knockout\_tourn(U[j])}$
    }
    \label{alg:gpsaf-end_beta}

    \vspace{1mm}
    \tcc{Replacement ($\gamma$)}
    \ForEach{$j\gets1$ \KwTo $\tt{size(P)}$}{
    \label{alg:gpsaf-start_rho}
    
            $\rho \gets \tt{repl\_prob}(\tt{U[j]}, U, \gamma)$
    
            \lIf{$\tt{rand}() < \rho$}{
                $\tt{P[j]} \gets V[j])$
            }

    }
     \label{alg:gpsaf-end_rho}

    \vspace{1mm}
    \tcc{Evaluate on ESE}
    $\tt{P.F, P.G} \gets \tt{evaluate(P.X)}$
    \label{alg:gpsaf-evaluate}
    
    \vspace{1mm}
    \tcc{Prepare next iteration of GPSAF}
    $\Phi\tt{.advance(P.X, P.F, P.G)}$
    \label{alg:gpsaf-advance}

    $\tt{e} \gets \tt{update\_error(P.F, P.G, P.}\hat{F}, \tt{P.}\hat{G})$

    \vspace{1mm}
    $\tt{A} \gets \tt{A} \cup \tt{P}$

%\vspace{1mm}
}

\caption{GPSAF: \textbf{G}eneralized \textbf{P}robablistic \textbf{S}urrogate-\textbf{A}ssisted \textbf{F}ramework}
\label{alg:gpsaf}
\end{algorithm}}

% general goals and ideas behind the method
However, how shall this interface now be utilized, and what role can surrogates play in improving the algorithm's performance? Precisely this is the subject of this article. Nevertheless, before moving on to the proposed framework, some more specifications of the surrogate usage are to be defined: First, the surrogate shall only be used as an assistant (in contrast to other methods where everything is developed centered around the surrogate). Second, the proposed method should be adaptive, allowing to decrease and increase the impact of surrogate usage and, if desired, even falling back to the original pseudo-code shown in Algorithm~\ref{alg:infill-advance}. Third, the surrogate prediction error needs to be addressed to ensure both exploitation and exploration. Altogether, the design goals are formulated to make the optimization framework and surrogate incorporation flexible.
The proposed Generalized Probabilistic Surrogate-Assisted Framework (GPSAF) meets these goals by introducing two different phases: First, the $\alpha$-phase using the current state of algorithm $\Phi$ introducing some surrogate influence. This pre-filtering phase uses a replacement strategy based on surrogate predictions;
second, the $\beta$-phase continues to run algorithm $\Phi$ for \emph{multiple} consecutive iterations on the surrogate resulting in a (convergence) search pattern. From the search pattern, solutions are selected by applying a probabilistic knockout tournament. The tournament incorporates the distribution of historical prediction errors by comparing solutions under noise to address the surrogate inaccuracies. The two phases and their control parameter allow configuring the surrogate usage, and the probabilistic knockout tournament serves as a self-adaptive mechanism to balance its exploitation and exploration.

\subsection{Generalized Probabilistic Surrogate-Assisted Framework}

% pseudo-code and how the method is explained in the following
Before describing the responsibilities and details of each of the phases, the outline of the algorithm shall be discussed (see Algorithm~\ref{alg:gpsaf}).
Before any surrogate can be fit, a solution archive $A$ is initialized by some design of experiments $\tt{A.X}$ are generated in a space-filling manner.
A good spread of solutions is recommended to allow surrogates to capture the overall fitness landscape as accurately as possible.
$\tt{A.X}$ is evaluated on the expensive solution evaluation (ESE) resulting in $\tt{A.F}$ and $\tt{A.G}$ (Line~\ref{alg:gpsaf-doe}).
Then, while the number of evaluations is less than the maximum solution evaluation budget $\tt{SE}^{(\max)}$, infill solutions $\tt{P.X}$ are generated by calling the non-deterministic \emph{infill} method of the baseline algorithm~$\Phi$. 
The default execution of algorithm $\Phi$ would immediately evaluate $\tt{P.X}$ using ESE and directly feed the solutions back to the algorithm by executing $\Phi\tt{.advance(P.X, P.F, P.G)}$ (Line~\ref{alg:gpsaf-evaluate} and~\ref{alg:gpsaf-advance}). However, instead of doing so, GPSAF modifies $P.X$ in a way to be influenced and biased by surrogates (Line~\ref{alg:gpsaf-fit-surr} to~\ref{alg:gpsaf-end_beta}) and advances the algorithm in the end of the iteration (Line~\ref{alg:gpsaf-evaluate} to~\ref{alg:gpsaf-advance}).
After having estimated the surrogate error and fitted the surrogates for objective and constraint functions, the $\alpha$-phase adds surrogate influence to $\tt{P.X}$  by replacing solutions being predicted to be better (Line~\ref{alg:gpsaf-start_alpha} to~\ref{alg:gpsaf-end_alpha}). 
Thereafter, the $\beta$-phase runs algorithm $\Phi$ for multiple generations (evaluations only on ASE) and assigns each solution to its closest $\tt{P.X}$. For each of the resulting candidate solution pools $\tt{U[j]}$ assigned $\tt{P[j]}$ a probabilistic tournament determines the winning candidate (Line~\ref{alg:gpsaf-start_beta} to~\ref{alg:gpsaf-end_beta}). 
Afterward, the replacement phase takes place where either the solution originating from the  $\alpha$-phase $\tt{P[j]}$ is kept or replaced with $\tt{U[j]}$ from the $\beta$-phase (Line~\ref{alg:gpsaf-start_rho} to~\ref{alg:gpsaf-end_rho}).
The solutions set to $\tt{P.X}$ are evaluated, and the algorithm $\Phi$ is advanced (Line~\ref{alg:gpsaf-evaluate} and~\ref{alg:gpsaf-advance}). Finally, the prediction error is updated before starting the next iteration, and the newly evaluated solutions are added to archive A.

% transition to next sub-sections
The overall outline of GPSAF shall provide an idea of where and when the $\alpha$ and $\beta$ phases take place and what role they play in modifying the infill solutions fed back to the algorithm. Next, each phase shall be explained and discussed in detail.

\subsection{Surrogate Influence through Tournament Pressure ($\alpha$)}
\label{sec:alpha}

The first mechanism of GPSAF incorporates tournament pressure by utilizing the predictions of a surrogate as a referee. Tournament pressure is a well-known concept in evolutionary computation to introduce a bias towards more promising solutions~\cite{1995-miller-tournament-selection}. 
Usually, its purpose is to introduce selection bias during mating to increase the chances of involvement of better-performing individuals. 
A specifically helpful control parameter is the number of competitors in each tournament to naturally increase or decrease the selection pressure.

\begin{figure}[t]
    \centering
    % trim=left bottom right top, clip. (add fbox)
    \includegraphics[page=2,width=0.80\linewidth,trim=0 3.5cm 19.5cm 0,clip]{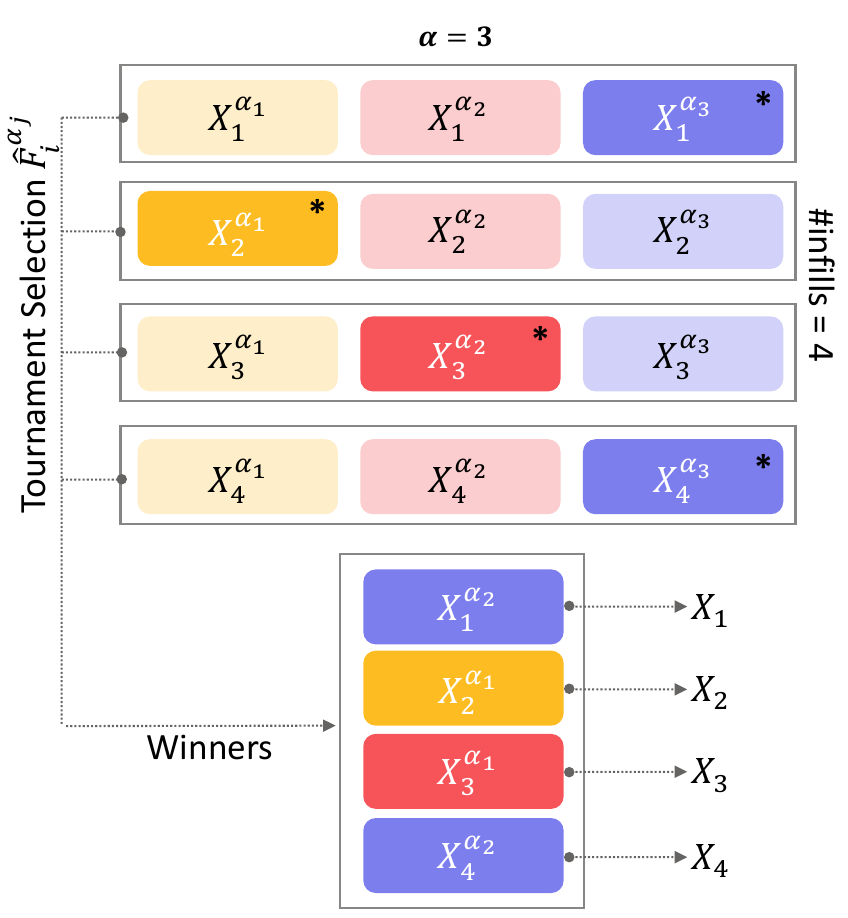}
    \caption{Tournament selection with $\alpha$ competitors to create a surrogate-influenced infill solutions.}
    \label{fig:alpha}
\end{figure}

Here, we borrow the tournament pressure mechanism to provide solutions with \emph{surrogate} influence before their evaluation on the EOP. 
The surrogate predictions provide the necessary ASEs to determine the winner of each tournament.
The number of competitors $\alpha$ in each tournament can control the amount of surrogate influence.
In Figure~\ref{fig:alpha} the surrogate-assisted tournament selection with three competitors ($\alpha=3$) for four infill solutions ($n=4$) is visualized. Initially, the algorithm's \emph{infill} function is called three times to generate the solution sets $X^{\alpha_1}$, $X^{\alpha_2}$, and $X^{\alpha_3}$.
After evaluating each of the solution sets on the surrogate, a tournament takes place where $\alpha$ solutions of the $j$-th infill solution set $X^{\alpha_j}$ compete with each other. For instance, for the first tournament, the winner of $X^{\alpha_1}_1$, $X^{\alpha_1}_2$, and $X^{\alpha_1}_3$ is declared. 
The winner of each solution pool is determined as follows: if \emph{all} solutions are infeasible, select the least infeasible solution; otherwise, select a non-dominated solution (break ties randomly).
For both the constraint and objective values, only ASEs are used. 
By repeating the tournament $n$ times and declaring the winners $X_i$ where $i \in (1,\ldots,n)$, here $X_1$, $X_2$, $X_3$, and $X_4$, four surrogate-influenced solutions have been selected and the $\alpha$-phase is completed.

% implementation and configuration
It is worth noting that because the \emph{infill} method calls of algorithm $\Phi$ is non-deterministic and do not have any order, this can be implemented memory-friendly by a for loop as shown in Algorithm~\ref{alg:gpsaf} (see Line~\ref{alg:gpsaf-start_alpha} to~\ref{alg:gpsaf-end_alpha}). 
Generally, it shall also become apparent that setting $\alpha=1$ disables the tournament selection and serves as a fallback to the original algorithm. By involving the surrogate in the tournament selection ($\alpha > 1$), the infill solutions $\tt{P}$ get a smaller or larger influence based on the number of competitors, which provides a natural inclusion of surrogate guidance.

\subsection{Continue Optimization on Surrogate ($\beta$)}

% why is the beta-phase necessary, and what is the idea
After completing the $\alpha$-phase, the solution set $\tt{P}$ is already influenced by surrogates. But what are the limitations of the $\alpha$-phase, and why is there a necessity for a second one?
Even though the \emph{infill} method is called multiple times, the algorithm is never \emph{advanced} to the next iteration. Thus, it provides a \emph{one step} look-ahead, which is not sufficient to find near-optimal solutions on the surrogate (which to some degree can be very useful for convergence).
So, to further increase the surrogate's impact, the second phase looks $\beta$ iterations into the future by calling infill \emph{and} advance of the baseline algorithm repetitively.
Whereas for a smaller $\beta$, the surrogates will be somewhat exploited, for a larger $\beta$, near-optimal solutions using ASEs will be found (similarly to optimizing the surrogate directly).
However, it must be considered that ASEs have an underlying prediction error and can not be taken for granted.
In GPSAF, the error is addressed by making not only use of the final solution set resulting from the $\beta$ optimization runs but using the whole \emph{search pattern}. 
The pattern is first divided into multiple clusters by assigning all solutions to their closest solution (in the design space) to each solution in $\tt{P}.X$.
Then, we use a so-called probabilistic knockout tournament (PKT) to select solutions from each cluster with the goal of self-adaptively exploiting surrogates. The goal is to use surrogates more when they provide accurate predictions but use them more carefully when they provide only rough estimations. 
Necessary for generalization, PKT also applies to problems with multiple objectives and constraints, often with varying complexities and surrogate errors to be considered.

{\small
% \linespread{1.0}\selectfont
\begin{algorithm}[!t]
\SetKwInput{Input}{Input~}
\Input{Solution Set $C$, Prediction errors $e$, Number of winners $k$
}

\vspace{1mm}

$C^{(1)} \gets \tt{shuffle}(C)$
\label{alg:gpsaf-pkt-shuffle}

$t \gets 1$

\While{$|C^{(t)}| > k$}{

\lIf{$|C^{(t)}|$ is \tt{odd}}{$C^{(t)} \gets C^{(t)} \, \cup \, \tt{rselect(C^{(t)}, 1)}$}
\label{alg:gpsaf-pkt-odd}

$C^{(t+1)} \gets \emptyset$

\ForEach{$i\gets1$ \KwTo $|C^{(t)}| / 2$}{

$w \gets \tt{compare\_noisy}(C^{(t)}_{2i}, C^{(t)}_{2i + 1}, e)$
\label{alg:gpsaf-pkt-cmp}

$C^{(t+1)} \gets C^{(t+1)} \cup w$
\label{alg:gpsaf-pkt-add}

}

$t \gets t + 1$

}

\If{$|C^{(t)}| < k$}{

$C^{(t)} \gets C^{(t)} \cup \tt{rselect}(C^{(t-1)} \setminus C^{(t)}, |C^{(t)}| - k)$
\label{alg:gpsaf-pkt-finalize}

}

\Return $C^{(t)}$

\caption{Probabilistic Knockout Tournament (PKT)}
\label{alg:pkt}
\end{algorithm}}

Generally, we define PKT as a subset selection of $k$ solutions from a set of solutions $C$ by applying \emph{pairwise} comparisons under noise as shown in Algorithm~\ref{alg:pkt}.
Initially, the solution set $C$ to select from is shuffled to randomize the matches (Line~\ref{alg:gpsaf-pkt-shuffle}). If the current number of participants $|C^{(t)}|$ is odd, a random solution is chosen to compete twice (Line~\ref{alg:gpsaf-pkt-odd}).
Each competition occurs under noise, based on the current prediction error of the surrogates. The noise is added to each objective and constraint independently before comparing the solutions. After adding the noise, the comparison is identical to the subset selection explained in Section~\ref{sec:alpha} (feasibility, dominance, random tie break) with two competitors ($\alpha=2$). The winner of each round moves on to the next and is added to $C^{(t+1)}$ (Line~\ref{alg:gpsaf-pkt-add}).
Finally, if too many solutions have been eliminated, randomly choose some losers from the last round (Line~\ref{alg:gpsaf-pkt-finalize}). This results in a set of solutions of size $k$ being returned as tournament winners under noise.
The design of PKT applies to the most general case of constrained multi-objective optimization because the selection procedure can be reduced two a comparison of two solutions.

\begin{figure}[t]
    \centering\includegraphics[page=3,width=0.83\linewidth,trim=0 0cm 15.5cm 0,clip]{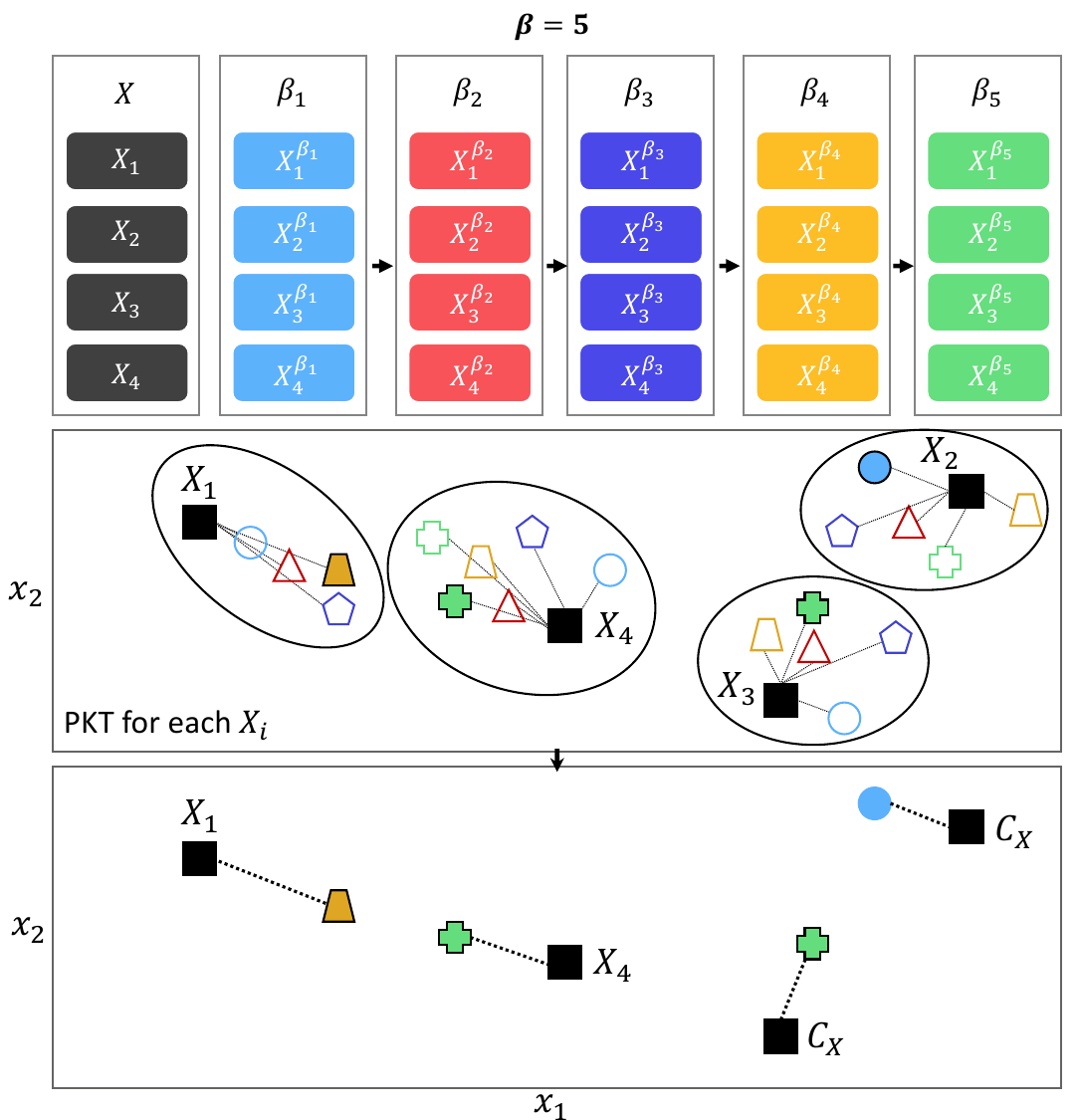}
    \caption{Continue Running the Algorithm for $\beta$ Iteration on the Surrogate.}
    \label{fig:beta}
\end{figure}

Back to the cluster-wise selection in the $\beta$-phase where PKT is executed with $k=1$ to obtain a winner for each solution set $\tt{U_j}$.
An example with five iterations ($\beta = 5$) and four infill solutions $X_1$, $X_2$, $X_3$, and $X_4$ is illustrated in Figure~\ref{fig:beta}. Calling the \emph{infill} and \emph{advance} function of the baseline algorithm results in five solution sets ($\beta_1$ to $\beta_5$) with four solutions each. The advancement of multiple iterations is based on ASEs.
In each iteration, all solutions are directly assigned to the closest $X_i$ solution from the $\alpha$-phase forming the cluster $U_i$. The cluster search pattern division is essential to preserve diversity.
For each cluster, a winner $V_i$ is declared by performing the PKT. For instance, in this example, $X_1$ has four solutions in $U_1$ where one from the fourth iteration $\beta_4$ is finally selected.
At the end of the $\beta$-phase, each cluster $U_i$ has at most one solution $V_i$ to be assigned to (some clusters may stay empty because no solutions are assigned to it).

The $\beta$-phase exploits the surrogates significantly more than the $\alpha$-phase by optimizing multiple iterations on the surrogates. In addition, mechanisms such as cluster-based search pattern selection help preserve diversity. Finally, it is worth noting that, analogously to the $\alpha$-phase, the surrogate assistance can be disabled by setting a specific configuration ($\beta=0$). 
Thus, GPSAF provides a fallback mechanism to the baseline optimization method without surrogate assistance (when $\alpha=1$ and $\beta=0$). Increasing one or the other will add more and more guidance through surrogate models.
How the two phases, each resulting in a set of solutions, are now combined to find a trade-off between the more explorative first and more exploiting second phase shall be discussed next.

\subsection{Balancing the Exploration and Exploitation ($\gamma$)}
\label{sec:rho}

The $\alpha$ and $\beta$ phases are designed to add surrogate assistance to an algorithm. The $\beta$ phase utilizes the search pattern on the surrogate and assigns solutions from the $\alpha$-phase. This means one can now choose to either stick to the more explorative $\alpha$ or exploit the $\beta$ solution. 
The most simple way of making this choice is by replacing the solution with probability $\rho$. 
However, in a pilot study addressing single-objective optimization, it has been shown that a more dynamic selection strategy is beneficial~\cite{2021-blankjul-psaf}. 

A piece of particularly useful information for making this decision is the distribution of assigned solutions across clusters.
The search pattern derived from surrogates with a high-density area indicates a region of interest. Thus, we propose to set the replacement probability:

\begin{equation}
    \rho = \left( \frac{|U_j|}{\max_j |U_j|} \right) ^ \gamma
\end{equation}

The denominator $\max_j |U_j|$ normalizes the number of assigned points with respect to the points in the current cluster $|U_j|$. The exponent $\gamma$ can be used to control the importance of the distribution and was kept constant at  $\gamma = 0.5$. 
The cluster with the highest density is always chosen from the $\beta$-phase because the nominator and denominator will be equal. 
This will necessarily be the case for baseline algorithms returning only one infill solution where a stronger surrogate bias is desirable. 
After the replacement, the solutions will finally be sent to the time-consuming solution evaluation.

\subsection{Surrogate Management}

Besides using surrogates in an algorithmic framework, some more words need to be said about the models themselves.  
First, one shall note that only the predictions of data points need to be provided by surrogates and no additional error estimation (the error estimates are kept track of by our method directly). Not requiring an error estimation does not limit the models to a specific type, unlike other surrogate-based algorithms.
Second, each of the objective and constraint functions is modeled independently, known as \emph{M1} in the surrogate usage taxonomy in~\cite{2019-deb-taxonomy}.
Even though modeling all functions increases the algorithmic overhead, it prevents larger prediction errors through complexity aggregations of multiple functions. 
Third, a generic framework for optimizing computationally expensive functions requires a generic surrogate model implementation. 
Clearly, some model types are more suitable for some problems than others. Thus, to provide a more robust framework, each function is approximated with a set of surrogates, and the best one is finally chosen to be used.
The surrogate types in this paper consist of the model types RBF~\cite{1971-hardy-rbf} and Kriging~\cite{1951-krige-mine}, both initialized with different hyper-parameters (normalization, regressions, kernel). 
A pre-normalization step referred to as PLOG~\cite{2017-bagheri-sacobra} is attempted and selected if well-performing for constraint functions.
Two metrics assess the performance of a model: First, Kendall Tau Distance~\cite{1938-kendall-tau} comparing the ranking of solutions being less sensitive to outliers with a large prediction error; second, the Maximum Absolute Error (MAE) to break any ties. The value of MAE is also used as an error approximation when noise is added to individuals.
The error estimation in the first iteration is based on k-fold cross-validation ($k=5$) to get a rough estimate of how well a surrogate can capture the function type. The performance metrics are updated in each iteration by taking all solutions seen so far as training and the newly evaluated solutions as a test set. Finally, a moving average of five iterations to avoid a smooth and more robust estimation provides the data for selecting the best surrogate and estimating the prediction error for each objective and constraint.

\section{Experimental Results}
\label{sec:results}

\input{tables/soo}

In this section, we present the performance of GPSAF applied to various population-based algorithms solving unconstrained and constrained, single- and multi-objective optimization problems.
Proposing an optimization framework requires comparing a group of algorithms, which is not a trivial task itself. Benchmarking is further complicated when non-deterministic algorithms are compared, in which case not only a single but multiple runs need to be considered.

For a fair comparison of optimization methods across test problems and to measure the impact of GPSAF on a baseline algorithm, we use the following ranking-based procedure:

\vspace{1mm}

\begin{enumerate}[label=\textbf{\roman*.}, itemsep=0pt]

\item \textbf{Statistical Domination:}  After collecting the data for each test problem and algorithm ($\mathcal{A} \in \Omega$) from multiple runs, we perform a pairwise comparison of performance indicators (PI) between all algorithms using the Wilcoxon Rank Sum Test ($\alpha=0.05$). The null-hypothesis $H_0$ is that no significant difference exists, whereas the alternative hypothesis is that the performance indicator of the first algorithm ($\texttt{PI}(\mathcal{A})$) is smaller than the one of the second one ($\texttt{PI}(\mathcal{B})$)
The PI is for single-objective optimization the gap to the optimum (if known) or the best function value found. For multi-objective optimization IGD~\cite{2004-coello-igd} (if optimum is known) or Hypervolume~\cite{1998-zitzler-hv} is used.
\begin{equation}
    \phi(\mathcal{A}, \mathcal{B}) = \tt{\textbf{RANKSUM}}(\tt{PI}(\mathcal{B}), \tt{PI}(\mathcal{A}),  \tt{alt='less'}),
\end{equation}

where the function $\phi(\mathcal{A}, \mathcal{B})$ returns zero if the null hypothesis is accepted or a one if it is rejected. 

\vspace{1mm}
\item \textbf{Number of Dominations}: The performance $P(\mathcal{A})$ of algorithm~$\mathcal{A}$ is then determined by the number of methods that are dominating it:
\begin{equation}
    P(\mathcal{A}) = \sum_{\substack{\mathcal{B} \in \Omega\\ \mathcal{A} \neq \mathcal{B}}} \phi(\mathcal{B}, \mathcal{A})
\end{equation}
This results in a domination number $P(\mathcal{A})$ for each method, which is zero if no other algorithm does not outperform it.

\vspace{1mm}
\item \textbf{Ranking}: Finally, we sort the methods by their $P(\mathcal{A})$. This may result in a partial ordering with multiple algorithms with the same $P(\mathcal{A})$ values. In order to keep the overall sum of ranks equal, we assign their average ranks in case of ties.
For instance, let us assume five optimizations methods A, B, C, D, and E:  algorithm A outperforms all others; between the performances of B, C, and D, no significant difference exists; E performs the worst. In this case, method A gets rank $1$, the group of methods B, C, and D, rank $(2 + 3 + 4) / 3 = 9/3 = 3$, and E rank $5$.  Averaging the ranks for ties penalizes an optimization method for being dominated by the same amount of algorithms as others and keeps the rank sum for each problem the same.

\end{enumerate}

This conveniently provides a ranking for each test problem. To evaluate the performance of a method on a test suite, we finally average the ranks across problems. 
If an algorithm fails to solve a specific problem for all runs, it gets the maximum rank and becomes the worst performing algorithm. Otherwise, all failing runs will be ignored (this has only rarely happened for a competitor algorithm to compare with).
The ranks shall be used to compare the performances of methods in this manuscript, the values of the performance indicators for the methods on all test problems can be found in the Supplementary Document.
Each algorithm has been executed 11 times on each test problem. If not explicitly mentioned in the specific experiment, the total number of solution evaluations has been set to $\tt{SE}^{(\max)}=300$. For some simpler constrained problems, even fewer evaluations have been used. 
A relatively limited evaluation budget also means that more complicated problems might not be solved (near) optimally.
However, a comparison of how well an algorithm has performed shall imitate the situation researchers face in practice.
If the number of variables is not fixed, the number of variables is fixed to $10$.
The results are presented in ranking tables where the overall best performing algorithm(s) are highlighted with a gray cell background for each ranking-based comparison for a test problem. The best-performing ones in a group are shown in bold.

Moreover, some more details about our implementation shall be said. For the baseline algorithms, we use implementations of population-based algorithms available in the well-known multi-objective optimization framework pymoo\footnote{http://pymoo.org (Version 0.5.0)}~\cite{2020-blank-pymoo} developed in Python. For all methods, the default parameters provided by the framework are kept unmodified, except the population size (=20) and the number of offsprings (=10) to create a more greedy implementation of the methods.
The surrogate implementation of Kriging is based on a Python clone\footnote{https://pypi.org/project/pydacefit/} of DACEFit~\cite{2002-lophaven-DACE} originally implemented in Matlab. The RBF models are a re-implementation based on~\cite{2017-bagheri-sacobra}.
The hyper-parameters of GPSAF were determined through numerous empirical experiments during the algorithm development.
A reasonable and well-performing configuration given by $\alpha=30$, $\beta=5$, and $\gamma=0.5$ is fixed throughout all experiments.

\subsection{(Unconstrained) Single-objective Optimization}

\input{tables/csoo}

The first experiment investigates the capabilities of GPSAF for improving the performance of existing algorithms on unconstrained single-objective problems. 
We use the BBOB test problems (24 functions in total) available in the COCO-platform~\cite{2020-hansen-coco} which is a widely used test suite with a variety of more and less complex problems.
Four well-known population-based optimization methods, DE~\cite{1997-storn-de}, GA~\cite{1989-goldberg-ga}, PSO~\cite{1995-kennedy-pso}, and CMAES~\cite{2001-hansen-cmaes} serve as baseline optimization algorithms and their GPSAF variants provide a surrogate-assisted version.
The results are compared with four other surrogate-assisted algorithms, SACOSO~\cite{2017-sun-SACOSO}, SACC-EAM-II~\cite{2019-blanchard-SACCEAMII}, SADESammon~\cite{2020-chen-SADESammon}, SAMSO~\cite{2021-fan-samso} available in the PlatEMO~\cite{2017-platemo} framework.
The rankings from the experiment are shown in Table~\ref{tbl:results-soo}. First, one can note that \mbox{GPSAF} outperforms the other four existing surrogate-assisted algorithms. One possible reason for the significant difference could be their development for a different type of test suite (for instance, problems with a larger number of variables). In this test suite, some problems are rather complicated, and exploiting the surrogate too much will cause to be easily trapped in local optima. Also, we contribute the efficiency of GPSAF to the significant effort for finding the most suitable surrogate. The order of relative rank improvement is given by GA ($6.292 / 2.292 = 2.7452$), PSO ($2.0927$), DE ($1.9546$), and for CMAES ($1.4522$). Besides GPSAF-GA having the biggest relative rank improvement, it also is the overall best performing algorithm in this experiment, closely followed by GPSAF-CMAES.
Altogether, a significant and quite remarkable improvement is achieved by applying GPSAF for (unconstrained) single-objective optimization.

\subsection{Constrained Single-objective Optimization}

Rarely are optimization problems unconstrained in practice. Thus, especially for surrogate-assisted methods aiming to solve computationally expensive real-world problems, the capability of dealing with constraints is essential.
The so-called G-problems or G-function benchmark~\cite{1996-michalewicz-g,1990-floudas-g} was proposed to develop optimization algorithms dealing with different kinds of constraints regarding the type (equality and inequality), amount, complexity, and result in feasible and infeasible search space. The original 13 test functions were extended in a CEC competition in 2006~\cite{2006-cec} to 24 constrained single-objective test problems~\cite{2006-liang-g}. In this study, G problems with only inequality constraints (and no equality constraints) are used.
Besides the GPSAF variants of DE and GA, improved stochastic ranking evolutionary strategy (ISRES)~\cite{2005-runarsson-isres} is applied to GSPAF. ISRES implements an improved mating strategy using differentials between solutions in contrast to its predecessor SRES~\cite{2000-runarsson-stochastic-ranking}.
ISRES follows the well-known 1/7 rule, which means with a population size of $\mu$ individuals $7 \cdot \mu$ offsprings are created. For this study, GPSAF creates a steady-state variant of ISRES by using the proposed probabilistic knockout tournament to choose \emph{one} out of the $\lambda$~solutions. 
This ensures a fair comparison with SACOBRA~\cite{2017-bagheri-sacobra} which also evaluates one solution per iteration. 
To the best of our knowledge, SACOBRA implemented in R~\cite{2020_R} is currently the best-performing algorithm on the G problem suite. 

The constrained single-objective results are presented in Table~\ref{tbl:results-csoo}. First, it is apparent that the GPSAF variants improve the baseline algorithms. Only for G2, the genetic algorithm outperforms its and other surrogate-assisted variants, which we contribute to the very restricted feasible search space (also, this has shown to be a difficult problem for surrogate-assisted algorithms in~\cite{2017-bagheri-sacobra}). 
Second, GPSAF-ISRES shows the best results out of all GPSAF variants. This indicates that it is beneficial if the baseline method has been proposed with a specific problem class in mind. Even though DE, GA, and PSO can handle constraints (for instance, naively using the parameter-less approach), there are known to not perform particularly well on complex constrained functions without any modifications. In contrast, ISRES has been tested on the G problems in the original study and proven to be effective. Furthermore, adding surrogate assistance to it has further improved the results.
Third, GPSAF-ISRES shows competitive performance to the state-of-the-art algorithm SACOBRA. In this experiment, out of all 14 test problems: GPSAF variants were able to outperform SACOBRA four times and a baseline algorithm (GA) one time; five times the performance of at least one GPSAF variant was similar; four times SACOBRA has shown significantly better results.
Altogether, one can say GPSAF has created surrogate-assisted methods competing with the state-of-the-art method for constrained single-objective problems.

\subsection{(Unconstrained) Multi-objective Optimization}

\input{tables/moo}
\input{tables/many}

Many applications have not one but multiple conflicting objectives to optimize. For this reason, this experiment focuses specifically on multi-objective optimization problems. As a test suite, we choose ZDT~\cite{2000-zitzler-zdt}, a well-known test suite proposed when multi-objective optimization has gained popularity. Throughout this experiment, we set the number of variables to $10$, except for the high multi-modal problem, ZDT4, where the number of variables is limited to $5$. The WFG~\cite{2005-wfg} test suit provides even more flexibility by being scalable with respect to the number of objectives. Here, we simply set the objective number to be two to create another bi-objective test suite. Moreover, the number of variables has been set to $10$ where four of them are positional.
The baseline algorithms NSGA-II~\cite{2002-deb-nsga2}, SMS-EMOA~\cite{2007-beume-sms-emoa}, and SPEA2~\cite{2001-zitzler-spea2} are used as baseline algorithms.
The results are compared with four other surrogate-assisted algorithms: AB-SAEA~\cite{2020-wang-absea}, KRVEA~\cite{2016-chugh-krvea}, ParEGO~\cite{2006-knowles-parego}, CSEA~\cite{2019-pan-csea} available in PlatEMO~\cite{2017-platemo}.

The results on the two multi-objective test suites are shown in Table~\ref{tbl:results-biobj}. First, one can note that all surrogate-assisted algorithms outperform the ones without. This indicates that surrogate assistance effectively improves the convergence behavior. 
Second, GPSAF-NSGA-II performs the best with a rank of $2.893$ and shows the best performance, followed by GPSAF-SPEA2, GPSAF-SMS-EMOA, and KRVEA. It is worth noting that ParEGO is penalized by being terminated for ZDT4 and WFG2, where the surrogate model was not able to be built.

To show the behavior of three-objective optimization problems, we have replaced NSGA-II with NSGA-III and run all algorithms on the DTLZ problems suite~\cite{2005-deb-dtlz} test suite. The results are shown in Table~\ref{tbl:results-many}. Whereas for most problems, the GPSAF variants outperform the baseline algorithms, for DTLZ1 and DTLZ3, this is not the case. Both problems consist of multi-modal convergence functions, which causes a large amount of surrogate error. Thus, surrogate-assisted algorithms (including the four GPSAF is compared to) are misguided. This seems to be a vital observation deserving to be investigated in more detail into the future. Nevertheless, GPSAF improves the performance of baseline algorithms for the other problems. GPSAF-SMS-EMOA shows overall the best results in this experiment with an average rank of $2.786$ followed by GPSAF-NSGA-III.

\subsection{Constrained Multi-objective Optimization}

\input{tables/cmoo}

Lastly, we shall compare GPSAF on constrained multi-objective optimization problems which often occur in real-world optimization. The challenge of dealing with multiple objectives and constraints in combination with computationally expensive solution evaluations truly mimics the complexity of industrial optimization problems.
We have compared our results with HSMEA~\cite{2019-habib-hsmea} a recently proposed algorithm for constraint multi-objective optimization. 
With consultation of the authors, some minor modifications of the publicly available source code had to be made for dealing with computationally expensive constraints -- as this is an assumption made in this study.
The results on CDTLZ~\cite{2014-deb-nsga3-part2}, BNH~\cite{1997-binh-bnh}, SRN~\cite{1994-srinivas-srn}, TNK~\cite{1995-tanaka-tnk}, and OSY~\cite{1995-osyczka-osy} are shown in~ Table~\ref{tbl:results-constr-multi}.
Again, one can observe that the GPSAF variants consistently improve the performance of the baseline optimization methods. The only exception is C1-DTLZ1, where all methods could find no feasible solution, and thus, an equal rank is assigned. We contribute this to the complexity of the test problems given by the constraint violation and the multi-modality of the objective functions.
For OSY, TNK, the GPSAF variants show a significantly better performance than HSMEA; for C3-DTLZ4, the performance is similar; and for C2-DTLZ2, BNH, and TNK, it performs better.
Altogether, GPSAF-NSGA-II can obtain a better rank than HSMEA, but it shall be fair to say that for three out of the seven constrained multi-objective optimization problems, HSMEA is the winner. 
Nevertheless, GPSAF improved the performance of baseline algorithms and showed competitive results to another surrogate-assisted optimization method.

\section{Concluding Remarks}
\label{sec:conclusion}

% summary
This article has proposed a generalized probabilistic surrogate-assisted framework applicable to any type of population-based algorithm. GSPAF incorporates two different phases to provide surrogate assistance, one considering using the current state of the baseline algorithm and the other looking at multiple iterations into the future. In contrast to other existing surrogate-assisted algorithms, the surrogate search is not reduced to the final solutions on the surrogate, but the whole search pattern is utilized. Solutions are selected using a probabilistic tournament that considers surrogate prediction errors for objectives and constraints from the search pattern. 
GPSAF has been applied to multiple well-known population-based algorithms proposed for unconstrained and constrained single and multi-objective optimization. We have provided comprehensive results on test problem suites indicating that GPSAF competes and outperforms existing surrogate-assisted methods. The combination of GPSAF creating well-performing surrogate-assisted algorithms with its \emph{simplicity} and \emph{broad} applicability is very promising.

% hyper-parameters
The encouraging results provide scope for further exploring generalized surrogate-assisted algorithms. One main challenge of a generalized approach is the recommendation of hyper-parameter configurations ($\alpha$, $\beta$, $\rho$, or $\gamma$). The parameters have been set through empirical experiments; however, through the broad applicability, different mechanisms of baseline algorithms on very different optimization problems make it difficult to draw generally valid conclusions. A more systemic and possibly resource-intensive study shall provide an idea of how different hyper-parameter settings impact the performance of different algorithms. In addition, experiments investigating the sensitivity shall be especially of interest.

% more evaluations - how to select points
The focus of this study was to explore different types of problems with multiple objectives and constraints. Thus the number of variables was kept relatively small as this is often the case for computationally expensive problems. Thus, even though the search space dimensions do not directly impact the idea proposed in this article, it shall be part of a future study of how surrogate assistance performs for large-scale problems. Moreover, the number of solution evaluations per run has been set to $300$, which allows using all solutions exhaustively for modeling without a large modeling overhead. However, more solution evaluations might be feasible for mediocre expensive optimization problems.

Nevertheless, this extensive explorative study on the use of surrogates in single and multi-objective optimization with and without constraints has indicated a viable new direction in congruence with existing emerging studies for a generic optimization methodology.

% different variable types
%GPSAF does not make any assumption about variables types. Therefore, even though this paper has focused on continuous variables, the proposed idea applies to discrete or even mixed variables given suitable surrogates and optimization methods (following the \emph{infill-and-advance} interface). 
%Thus, a follow-up study investigating different types of surrogates designed for different types of variables shall be worth exploring.

%{\color{red} 
%The conclusion is more or less a draft. Feel free to modify it entirely. If we keep the constraint results, we might also mention that for unconstrained and constraint, multi and many more experiments are needed, and concepts such as a local search or eps constraint handling might help.
%}

% if have a single appendix:
%\appendix[Proof of the Zonklar Equations]
% or
%\appendix  % for no appendix heading
% do not use \section anymore after \appendix, only \section*
% is possibly needed

% use appendices with more than one appendix
% then use \section to start each appendix
% you must declare a \section before using any
% \subsection or using \label (\appendices by itself
% starts a section numbered zero.)
%

% use section* for acknowledgment
% \section*{Acknowledgment}

% Can use something like this to put references on a page
% by themselves when using endfloat and the captionsoff option.
\ifCLASSOPTIONcaptionsoff
  \newpage
\fi

% trigger a \newpage just before the given reference
% number - used to balance the columns on the last page
% adjust value as needed - may need to be readjusted if
% the document is modified later
%\IEEEtriggeratref{8}
% The "triggered" command can be changed if desired:
%\IEEEtriggercmd{\enlargethispage{-5in}}

% references section

% can use a bibliography generated by BibTeX as a .bbl file
% BibTeX documentation can be easily obtained at:
% http://mirror.ctan.org/biblio/bibtex/contrib/doc/
% The IEEEtran BibTeX style support page is at:
% http://www.michaelshell.org/tex/ieeetran/bibtex/
\bibliographystyle{IEEEtran}
\bibliography{IEEEabrv,references}
%
% <OR> manually copy in the resultant .bbl file
% set second argument of \begin to the number of references
% (used to reserve space for the reference number labels box)
% \begin{thebibliography}{1}

% \bibitem{IEEEhowto:kopka}
% H.~Kopka and P.~W. Daly, \emph{A Guide to \LaTeX}, 3rd~ed.\hskip 1em plus
%   0.5em minus 0.4em\relax Harlow, England: Addison-Wesley, 1999.

% \end{thebibliography}

% biography section
%
% If you have an EPS/PDF photo (graphicx package needed) extra braces are
% needed around the contents of the optional argument to biography to prevent
% the LaTeX parser from getting confused when it sees the complicated
% \includegraphics command within an optional argument. (You could create
% your own custom macro containing the \includegraphics command to make things
% simpler here.)
%\begin{IEEEbiography}[{\includegraphics[width=1in,height=1.25in,clip,keepaspectratio]{mshell}}]{Michael Shell}
% or if you just want to reserve a space for a photo:

% trim=left bottom right top, clip. (add fbox)
\begin{IEEEbiography}[{\includegraphics[width=1in,height=1.25in,clip,trim=2.5cm 15cm 2.5cm 0,keepaspectratio ]{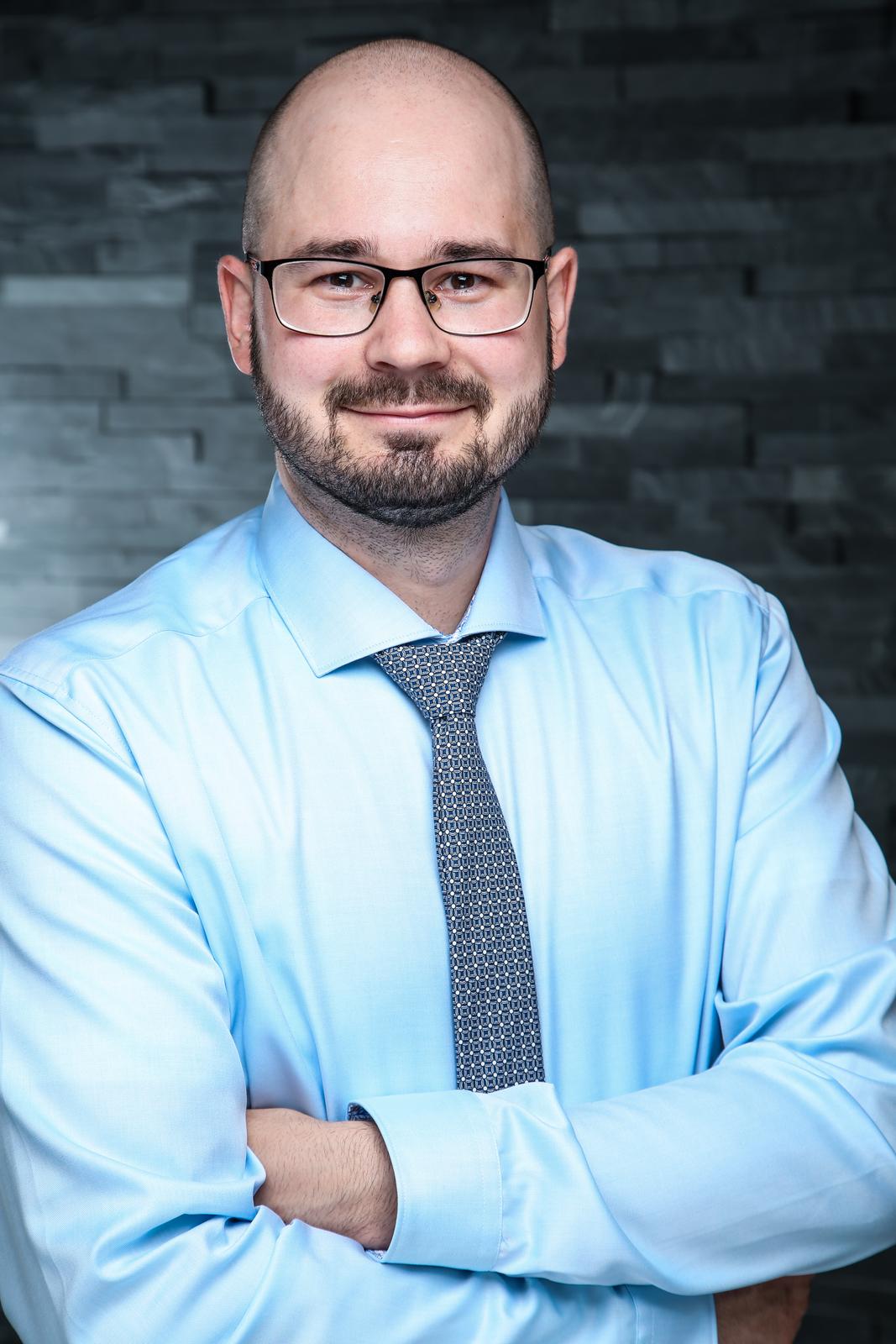}}]{Julian Blank}
is a PhD candidate in the Department of Computer Science and Engineering at Michigan State University. He received his B.Sc. in Business Information Systems from Otto von Guericke University, Germany, in 2010. He was a visiting scholar for six months at the Michigan State University, Michigan, USA, in 2015, and finished his M.Sc. in Computer Science at Otto von Guericke University, Germany, in 2016. He is the leading developer of pymoo, an open-source multi-objective optimization framework in Python. His research interests include evolutionary computation, multi-objective optimization, surrogate-assisted optimization, and machine learning.
\end{IEEEbiography}

% if you will not have a photo at all:
\begin{IEEEbiography}[{\includegraphics[width=1in,height=1.25in,clip,keepaspectratio]{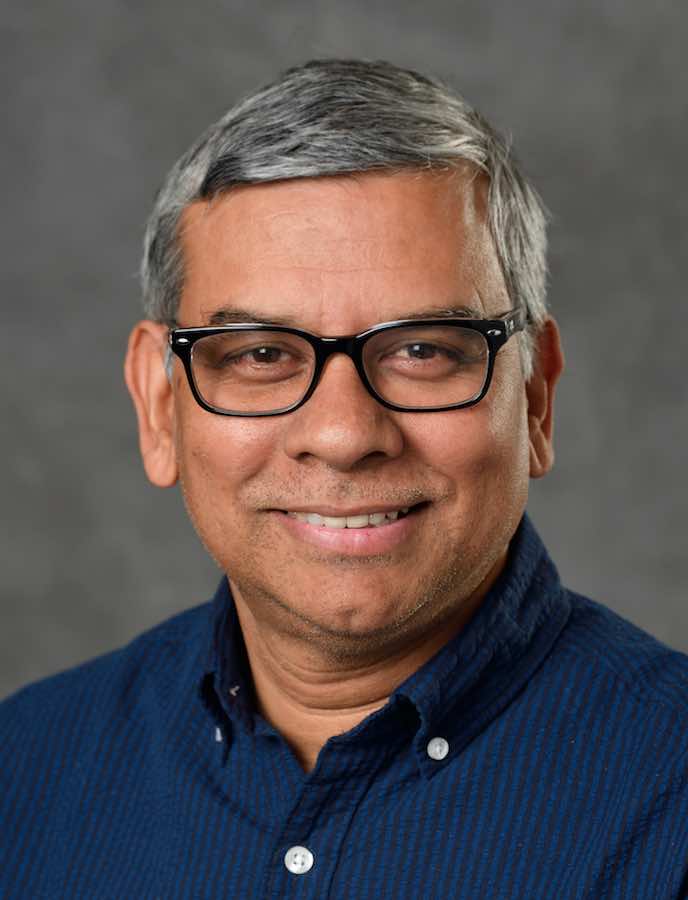}}]{Kalyanmoy Deb}
is Fellow, IEEE, University Distinguished Professor and Koenig Endowed Chair Professor with the Department of Electrical and Computer Engineering, Michigan State University, East Lansing, Michigan, USA. He received his Bachelor's degree in Mechanical Engineering from IIT Kharagpur in India, and his Master's and Ph.D. degrees from the University of Alabama, Tuscaloosa, USA, in 1989 and 1991. He is largely known for his seminal research in evolutionary multi-criterion optimization. He has published over 580 international journal and conference research papers to date. His current research interests include evolutionary optimization and its application in design, modeling, AI, and machine learning. He is one of the top-cited EC researchers with more than 190,000 Google Scholar citations.
\end{IEEEbiography}

% insert where needed to balance the two columns on the last page with
% biographies
%\newpage

% \begin{IEEEbiographynophoto}{Jane Doe}
% Biography text here.
% \end{IEEEbiographynophoto}

% You can push biographies down or up by placing
% a \vfill before or after them. The appropriate
% use of \vfill depends on what kind of text is
% on the last page and whether or not the columns
% are being equalized.

%\vfill

% Can be used to pull up biographies so that the bottom of the last one
% is flush with the other column.
%\enlargethispage{-5in}

% that's all folks
\end{document}

%% file: tables/soo.tex
\begin{table*}[!htb]
\caption{A comparison of DE, GA, PSO, and CMAES with their GPSAF variants on unconstrained single-objective problems with four other surrogate-assisted algorithms. 
The rank of the best performing algorithm in each group is shown in bold. The overall best performing algorithm for each problem is highlighted with a gray shade.
}
\centering
\begin{tabular}{|c|cc|cc|cc|cc|cccc|}
\toprule
Problem    &  DE &\specialcell{GPSAF-\\DE} &   GA  &  \specialcell{GPSAF-\\GA} & PSO & \specialcell{GPSAF-\\PSO} & CMAES & \specialcell{GPSAF-\\CMAES} & SACOSO & \specialcell{SACC-\\EAM-II } & SADESammon & SAMSO  \\
\midrule
f01 & 11.0 & \textbf{4.0} & 7.5 & \textbf{2.5} & 5.5 & \cellcolor{gray!25} \textbf{1.0} & 5.5 & \textbf{2.5} & \textbf{7.5} & 9.5 & 12.0 & 9.5 \\
f02 & 10.0 & \textbf{4.0} & 4.0 & \cellcolor{gray!25} \textbf{1.0} & 4.0 & \textbf{2.0} & 7.5 & \textbf{6.0} & 12.0 & \textbf{7.5} & 10.0 & 10.0 \\
f03 & 8.5 & \textbf{3.0} & 5.0 & \cellcolor{gray!25} \textbf{1.5} & 7.0 & \cellcolor{gray!25} \textbf{1.5} & 8.5 & \textbf{5.0} & 10.0 & \textbf{5.0} & 11.5 & 11.5 \\
f04 & 8.5 & \textbf{4.0} & 4.0 & \textbf{2.0} & 4.0 & \cellcolor{gray!25} \textbf{1.0} & \textbf{6.5} & \textbf{6.5} & 10.0 & \textbf{8.5} & 11.0 & 12.0 \\
f05 & 5.5 & \cellcolor{gray!25} \textbf{1.0} & 7.5 & \textbf{3.0} & 5.5 & \textbf{2.0} & 7.5 & \textbf{4.0} & 11.0 & \textbf{9.5} & \textbf{9.5} & 12.0 \\
f06 & \textbf{7.0} & \textbf{7.0} & \cellcolor{gray!25} \textbf{2.0} & \cellcolor{gray!25} \textbf{2.0} & \textbf{4.5} & \textbf{4.5} & 9.0 & \cellcolor{gray!25} \textbf{2.0} & 10.0 & \textbf{7.0} & 11.0 & 12.0 \\
f07 & 9.5 & \textbf{5.5} & 5.5 & \cellcolor{gray!25} \textbf{2.0} & 5.5 & \cellcolor{gray!25} \textbf{2.0} & 5.5 & \cellcolor{gray!25} \textbf{2.0} & 9.5 & \textbf{8.0} & 11.0 & 12.0 \\
f08 & 8.0 & \textbf{6.0} & 8.0 & \cellcolor{gray!25} \textbf{2.0} & 5.0 & \cellcolor{gray!25} \textbf{2.0} & 4.0 & \cellcolor{gray!25} \textbf{2.0} & 10.0 & \textbf{8.0} & 11.0 & 12.0 \\
f09 & 10.0 & \cellcolor{gray!25} \textbf{3.5} & 8.0 & \cellcolor{gray!25} \textbf{3.5} & \cellcolor{gray!25} \textbf{3.5} & \cellcolor{gray!25} \textbf{3.5} & 7.0 & \cellcolor{gray!25} \textbf{3.5} & \cellcolor{gray!25} \textbf{3.5} & 9.0 & 11.5 & 11.5 \\
f10 & 10.5 & \textbf{5.5} & 5.5 & \cellcolor{gray!25} \textbf{1.5} & 5.5 & \cellcolor{gray!25} \textbf{1.5} & \textbf{5.5} & \textbf{5.5} & 10.5 & \textbf{5.5} & 10.5 & 10.5 \\
f11 & 10.5 & \textbf{3.5} & 7.0 & \cellcolor{gray!25} \textbf{1.5} & 7.0 & \textbf{3.5} & \textbf{7.0} & \textbf{7.0} & 10.5 & \cellcolor{gray!25} \textbf{1.5} & 12.0 & 7.0 \\
f12 & \cellcolor{gray!25} \textbf{2.5} & 6.5 & \textbf{6.5} & \textbf{6.5} & \cellcolor{gray!25} \textbf{2.5} & \cellcolor{gray!25} \textbf{2.5} & \cellcolor{gray!25} \textbf{2.5} & 6.5 & 10.0 & \textbf{9.0} & 11.0 & 12.0 \\
f13 & 7.5 & \textbf{5.0} & 7.5 & \textbf{2.5} & 5.0 & \cellcolor{gray!25} \textbf{1.0} & 5.0 & \textbf{2.5} & 10.0 & \textbf{9.0} & 11.0 & 12.0 \\
f14 & 9.5 & \textbf{5.0} & 7.5 & \cellcolor{gray!25} \textbf{2.0} & 5.0 & \cellcolor{gray!25} \textbf{2.0} & 5.0 & \cellcolor{gray!25} \textbf{2.0} & 9.5 & \textbf{7.5} & 11.0 & 12.0 \\
f15 & 10.0 & \textbf{3.5} & 6.5 & \cellcolor{gray!25} \textbf{1.5} & 6.5 & \cellcolor{gray!25} \textbf{1.5} & 6.5 & \textbf{3.5} & 9.0 & \textbf{6.5} & 11.0 & 12.0 \\
f16 & 9.0 & \textbf{6.0} & 6.0 & \cellcolor{gray!25} \textbf{2.0} & 9.0 & \cellcolor{gray!25} \textbf{2.0} & 4.0 & \cellcolor{gray!25} \textbf{2.0} & 11.5 & \textbf{6.0} & 11.5 & 9.0 \\
f17 & 9.5 & \textbf{5.0} & 7.0 & \textbf{3.0} & 7.0 & \textbf{3.0} & 3.0 & \cellcolor{gray!25} \textbf{1.0} & 11.0 & \textbf{7.0} & 9.5 & 12.0 \\
f18 & 9.5 & \textbf{4.0} & 6.0 & \cellcolor{gray!25} \textbf{2.0} & 6.0 & \cellcolor{gray!25} \textbf{2.0} & 6.0 & \cellcolor{gray!25} \textbf{2.0} & 9.5 & \textbf{8.0} & 11.0 & 12.0 \\
f19 & 11.0 & \textbf{5.0} & 10.0 & \cellcolor{gray!25} \textbf{1.0} & \textbf{5.0} & \textbf{5.0} & 8.5 & \textbf{5.0} & 5.0 & \textbf{2.0} & 8.5 & 12.0 \\
f20 & 9.5 & \cellcolor{gray!25} \textbf{2.5} & 7.0 & \cellcolor{gray!25} \textbf{2.5} & \cellcolor{gray!25} \textbf{2.5} & \cellcolor{gray!25} \textbf{2.5} & \textbf{5.5} & \textbf{5.5} & 9.5 & \textbf{8.0} & 11.5 & 11.5 \\
f21 & 8.5 & \textbf{7.0} & \cellcolor{gray!25} \textbf{3.5} & \cellcolor{gray!25} \textbf{3.5} & \cellcolor{gray!25} \textbf{3.5} & \cellcolor{gray!25} \textbf{3.5} & \cellcolor{gray!25} \textbf{3.5} & \cellcolor{gray!25} \textbf{3.5} & 10.0 & \textbf{8.5} & 11.5 & 11.5 \\
f22 & 9.5 & \textbf{6.0} & 6.0 & \cellcolor{gray!25} \textbf{2.5} & 6.0 & \cellcolor{gray!25} \textbf{2.5} & \cellcolor{gray!25} \textbf{2.5} & \cellcolor{gray!25} \textbf{2.5} & 9.5 & \textbf{8.0} & 11.0 & 12.0 \\
f23 & 10.0 & \textbf{5.5} & 5.5 & \cellcolor{gray!25} \textbf{1.5} & \textbf{5.5} & \textbf{5.5} & \textbf{11.0} & 12.0 & 9.0 & \cellcolor{gray!25} \textbf{1.5} & 5.5 & 5.5 \\
f24 & 10.0 & \cellcolor{gray!25} \textbf{2.0} & 8.0 & \cellcolor{gray!25} \textbf{2.0} & 4.0 & \cellcolor{gray!25} \textbf{2.0} & 8.0 & \textbf{5.5} & 8.0 & \textbf{5.5} & 11.0 & 12.0 \\
\bottomrule
Total & 8.958 & \textbf{4.583} & 6.292 & \cellcolor{gray!25} \textbf{2.292} & 5.188 & \textbf{2.479} & 6.021 & \textbf{4.146} & 9.417 & \textbf{6.896} & 10.667 & 11.062 \\
\bottomrule
\end{tabular}
\label{tbl:results-soo}
\end{table*}

%% file: tables/csoo.tex
\begin{table*}[!htb]
\caption{A comparison of DE, GA, PSO, and ISRES with their GPSAF variants, on constrained single-objective problems with SACOBRA -- the current state-of-art algorithms for constrained optimization. 
}
\centering
\begin{tabular}{|cc|cc|cc|cc|c|c|}
\toprule
Problem    & $\tt{SE}^{(\max)}$ &  DE &\specialcell{GPSAF-\\DE} &   GA  &  \specialcell{GPSAF-\\GA} & PSO & \specialcell{GPSAF-\\PSO} & \specialcell{GPSAF-\\ISRES} & SACOBRA  \\
\midrule
G1 & 75 & 5.5 & \textbf{3.5} & 7.5 & \textbf{3.5} & 7.5 & \textbf{5.5} & \cellcolor{gray!25} \textbf{1.0} & \textbf{2.0} \\
G2 & 300 & \textbf{3.5} & 7.0 & \cellcolor{gray!25} \textbf{1.0} & 3.5 & 6.0 & \textbf{3.5} & \textbf{8.0} & \textbf{3.5} \\
G4 & 75 & 6.5 & \textbf{3.5} & 6.5 & \textbf{5.0} & 8.0 & \textbf{3.5} & \cellcolor{gray!25} \textbf{1.5} & \cellcolor{gray!25} \textbf{1.5} \\
G6 & 75 & 7.0 & \textbf{4.0} & 7.0 & \textbf{4.0} & 7.0 & \textbf{4.0} & \cellcolor{gray!25} \textbf{1.5} & \cellcolor{gray!25} \textbf{1.5} \\
G7 & 75 & 7.0 & \textbf{4.0} & 7.0 & \textbf{4.0} & 7.0 & \textbf{4.0} & \textbf{2.0} & \cellcolor{gray!25} \textbf{1.0} \\
G8 & 100 & 7.0 & \textbf{4.5} & 7.0 & \textbf{4.5} & 7.0 & \textbf{3.0} & \cellcolor{gray!25} \textbf{1.0} & \textbf{2.0} \\
G9 & 300 & 7.0 & \textbf{4.5} & 7.0 & \textbf{2.5} & 4.5 & \textbf{2.5} & \textbf{7.0} & \cellcolor{gray!25} \textbf{1.0} \\
G10 & 300 & 8.0 & \textbf{3.5} & 6.5 & \textbf{3.5} & 6.5 & \textbf{3.5} & \textbf{3.5} & \cellcolor{gray!25} \textbf{1.0} \\
G11 & 300 & 7.0 & \cellcolor{gray!25} \textbf{2.5} & 5.5 & \cellcolor{gray!25} \textbf{2.5} & 5.5 & \cellcolor{gray!25} \textbf{2.5} & \cellcolor{gray!25} \textbf{2.5} & \textbf{8.0} \\
G12 & 300 & 6.0 & \textbf{4.5} & 8.0 & \textbf{4.5} & 7.0 & \textbf{3.0} & \cellcolor{gray!25} \textbf{1.5} & \cellcolor{gray!25} \textbf{1.5} \\
G16 & 300 & 5.5 & \cellcolor{gray!25} \textbf{2.5} & \textbf{5.5} & 8.0 & 7.0 & \cellcolor{gray!25} \textbf{2.5} & \cellcolor{gray!25} \textbf{2.5} & \cellcolor{gray!25} \textbf{2.5} \\
G18 & 300 & 7.0 & \textbf{3.5} & 7.0 & \textbf{3.5} & 7.0 & \textbf{3.5} & \textbf{3.5} & \cellcolor{gray!25} \textbf{1.0} \\
G19 & 300 & 7.5 & \cellcolor{gray!25} \textbf{2.0} & 6.0 & \cellcolor{gray!25} \textbf{2.0} & 7.5 & \cellcolor{gray!25} \textbf{2.0} & \textbf{5.0} & \textbf{4.0} \\
G24 & 300 & 7.5 & \textbf{4.5} & 7.5 & \textbf{4.5} & 6.0 & \cellcolor{gray!25} \textbf{2.0} & \cellcolor{gray!25} \textbf{2.0} & \cellcolor{gray!25} \textbf{2.0} \\
\bottomrule
\multicolumn{2}{l}{Total} & 6.571 & \textbf{3.857} & 6.357 & \textbf{3.964} & 6.679 & \textbf{3.214} & \textbf{3.036} & \cellcolor{gray!25} \textbf{2.321} \\
\bottomrule
\end{tabular}
\label{tbl:results-csoo}
\end{table*}

%% file: tables/moo.tex
\begin{table*}[!htb]
\caption{A comparison of NSGA-II, SMS-EMOA, and SPEA2 with their GPSAF variants with four surrogate-assisted algorithms on bi-objective optimization problems.}
\centering
\begin{tabular}{|l|cc|cc|cc|cccc|}
\hline			
\toprule
Problem    &   NSGA-II  & \specialcell{GPSAF-\\NSGA-II} & SMS-EMOA &   \specialcell{GPSAF-\\SMS-EMOA}  & SPEA2 & \specialcell{GPSAF-\\SPEA2} & AB-SAEA & K-RVEA & ParEGO & CSEA \\
\midrule
ZDT1 & 9.0 & \cellcolor{gray!25} \textbf{2.5} & 9.0 & \cellcolor{gray!25} \textbf{2.5} & 9.0 & \cellcolor{gray!25} \textbf{2.5} & 6.0 & 5.0 & \cellcolor{gray!25} \textbf{2.5} & 7.0 \\
ZDT2 & 9.0 & \cellcolor{gray!25} \textbf{1.5} & 9.0 & \textbf{6.0} & 9.0 & \textbf{5.0} & 3.0 & 4.0 & \cellcolor{gray!25} \textbf{1.5} & 7.0 \\
ZDT3 & 9.0 & \textbf{4.5} & 9.0 & \textbf{4.5} & 9.0 & \textbf{4.5} & 4.5 & \cellcolor{gray!25} \textbf{1.0} & 2.0 & 7.0 \\
ZDT4 & 6.5 & \cellcolor{gray!25} \textbf{2.0} & 6.5 & \cellcolor{gray!25} \textbf{2.0} & 6.5 & \cellcolor{gray!25} \textbf{2.0} & 9.0 & 6.5 & 10.0 & \textbf{4.0} \\
ZDT6 & 7.5 & \textbf{3.0} & 9.5 & \textbf{6.0} & 9.5 & \textbf{3.0} & 5.0 & 3.0 & \cellcolor{gray!25} \textbf{1.0} & 7.5 \\
WFG1 & 9.0 & \textbf{6.0} & 9.0 & \textbf{6.0} & 9.0 & \textbf{6.0} & 4.0 & \cellcolor{gray!25} \textbf{2.0} & \cellcolor{gray!25} \textbf{2.0} & \cellcolor{gray!25} \textbf{2.0} \\
WFG2 & 7.0 & \cellcolor{gray!25} \textbf{1.5} & \textbf{8.0} & 9.0 & \textbf{4.5} & \textbf{4.5} & 4.5 & \cellcolor{gray!25} \textbf{1.5} & 10.0 & 4.5 \\
WFG3 & 8.0 & \textbf{3.5} & 10.0 & \textbf{3.5} & 8.0 & \cellcolor{gray!25} \textbf{1.5} & 5.5 & \cellcolor{gray!25} \textbf{1.5} & 5.5 & 8.0 \\
WFG4 & 6.5 & \cellcolor{gray!25} \textbf{1.5} & 9.0 & \textbf{5.0} & 6.5 & \cellcolor{gray!25} \textbf{1.5} & \textbf{3.5} & \textbf{3.5} & 9.0 & 9.0 \\
WFG5 & 8.0 & \textbf{2.5} & 8.0 & \textbf{4.0} & 10.0 & \textbf{5.0} & 6.0 & 2.5 & \cellcolor{gray!25} \textbf{1.0} & 8.0 \\
WFG6 & 9.0 & \cellcolor{gray!25} \textbf{2.0} & 10.0 & \textbf{6.5} & 4.5 & \cellcolor{gray!25} \textbf{2.0} & \cellcolor{gray!25} \textbf{2.0} & 4.5 & 6.5 & 8.0 \\
WFG7 & 5.5 & \textbf{4.0} & 8.5 & \cellcolor{gray!25} \textbf{2.0} & 8.5 & \cellcolor{gray!25} \textbf{2.0} & 5.5 & 7.0 & \cellcolor{gray!25} \textbf{2.0} & 10.0 \\
WFG8 & 7.5 & \cellcolor{gray!25} \textbf{2.5} & 10.0 & \textbf{5.0} & 9.0 & \cellcolor{gray!25} \textbf{2.5} & \cellcolor{gray!25} \textbf{2.5} & 6.0 & \cellcolor{gray!25} \textbf{2.5} & 7.5 \\
WFG9 & 7.5 & \cellcolor{gray!25} \textbf{3.0} & 7.5 & \cellcolor{gray!25} \textbf{3.0} & 6.0 & \cellcolor{gray!25} \textbf{3.0} & \cellcolor{gray!25} \textbf{3.0} & 10.0 & 9.0 & \cellcolor{gray!25} \textbf{3.0} \\
\bottomrule
Total & 7.786 & \cellcolor{gray!25} \textbf{2.857} & 8.786 & \textbf{4.643} & 7.786 & \textbf{3.214} & 4.571 & \textbf{4.143} & 4.607 & 6.607 \\
\bottomrule
\end{tabular}
\label{tbl:results-biobj}
\end{table*}

%% file: tables/many.tex
\begin{table*}[!htb]
\caption{A comparison of NSGA-III, SMS-EMOA, and SPEA2 with their GPSAF variants with four surrogate-assisted algorithms on three-objective optimization problems.}
\centering
\begin{tabular}{|l|cc|cc|cc|cccc|}
\hline			
\toprule
Problem    &   NSGA-III  & \specialcell{GPSAF-\\NSGA-III} & SMS-EMOA &   \specialcell{GPSAF-\\SMS-EMOA}  & SPEA2 & \specialcell{GPSAF-\\SPEA2} & AB-SAEA & K-RVEA & ParEGO & CSEA \\
\midrule
DTLZ1 & \cellcolor{gray!25} \textbf{1.5} & 3.5 & \cellcolor{gray!25} \textbf{1.5} & 6.5 & \textbf{5.0} & 8.0 & 10.0 & 9.0 & 6.5 & \textbf{3.5} \\
DTLZ2 & 8.5 & \cellcolor{gray!25} \textbf{2.0} & 8.5 & \cellcolor{gray!25} \textbf{2.0} & 8.5 & \cellcolor{gray!25} \textbf{2.0} & 6.0 & \textbf{4.0} & 5.0 & 8.5 \\
DTLZ3 & \cellcolor{gray!25} \textbf{2.0} & 5.5 & \cellcolor{gray!25} \textbf{2.0} & \cellcolor{gray!25} \textbf{2.0} & \textbf{5.5} & 8.0 & 10.0 & 9.0 & \textbf{5.5} & \textbf{5.5} \\
DTLZ4 & \textbf{6.5} & \textbf{6.5} & 9.0 & \textbf{6.5} & 6.5 & \textbf{3.5} & 2.0 & \cellcolor{gray!25} \textbf{1.0} & 10.0 & 3.5 \\
DTLZ5 & 6.0 & \textbf{3.0} & 9.0 & \cellcolor{gray!25} \textbf{1.0} & 9.0 & \textbf{2.0} & 6.0 & 6.0 & \textbf{4.0} & 9.0 \\
DTLZ6 & 7.5 & \textbf{5.5} & 7.5 & \textbf{3.5} & 9.0 & \textbf{5.5} & \cellcolor{gray!25} \textbf{1.5} & \cellcolor{gray!25} \textbf{1.5} & 10.0 & 3.5 \\
DTLZ7 & 8.0 & \textbf{3.5} & 8.0 & \textbf{3.5} & 8.0 & \textbf{3.5} & 3.5 & \cellcolor{gray!25} \textbf{1.0} & 10.0 & 6.0 \\
\bottomrule
Total & 5.714 & \textbf{4.214} & 6.5 & \cellcolor{gray!25} \textbf{3.571} & 7.357 & \textbf{4.643} & 5.571 & \textbf{4.5} & 7.286 & 5.643 \\
\bottomrule
\end{tabular}
\label{tbl:results-many}
\end{table*}

%% file: tables/cmoo.tex
\begin{table*}[!htb]
\caption{A comparison of NSGA-III, SMS-EMOA, and SPEA2 with their GPSAF variants with four surrogate-assisted algorithms on constrained multi-objective optimization problems.}
\centering
\begin{tabular}{|lc|cc|cc|cc|c|}
\toprule
Problem     & $\tt{SE}^{(\max)}$ &  NSGA-II  & \specialcell{GPSAF-\\NSGA-II} & SMS-EMOA &   \specialcell{GPSAF-\\SMS-EMOA}  & SPEA2 & \specialcell{GPSAF-\\SPEA2} & HSMEA\\
\midrule
C1-DTLZ1 & 300 & \cellcolor{gray!25} \textbf{4.0} & \cellcolor{gray!25} \textbf{4.0} & \cellcolor{gray!25} \textbf{4.0} & \cellcolor{gray!25} \textbf{4.0} & \cellcolor{gray!25} \textbf{4.0} & \cellcolor{gray!25} \textbf{4.0} & \cellcolor{gray!25} \textbf{4.0} \\
C2-DTLZ2 & 300 & \textbf{4.5} & \textbf{4.5} & \textbf{4.5} & \textbf{4.5} & \textbf{4.5} & \textbf{4.5} & \cellcolor{gray!25} \textbf{1.0} \\
C3-DTLZ4 & 300 & 5.0 & \cellcolor{gray!25} \textbf{1.5} & \textbf{5.0} & \textbf{5.0} & \textbf{5.0} & \textbf{5.0} & \cellcolor{gray!25} \textbf{1.5} \\
BNH & 100 & 5.5 & \textbf{2.5} & 7.0 & \textbf{4.0} & 5.5 & \textbf{2.5} & \cellcolor{gray!25} \textbf{1.0} \\
SRN & 100 & 6.0 & \textbf{3.0} & 6.0 & \textbf{3.0} & 6.0 & \textbf{3.0} & \cellcolor{gray!25} \textbf{1.0} \\
TNK & 100 & 6.0 & \cellcolor{gray!25} \textbf{2.0} & 6.0 & \cellcolor{gray!25} \textbf{2.0} & 6.0 & \cellcolor{gray!25} \textbf{2.0} & \textbf{4.0} \\
OSY & 300 & 5.0 & \cellcolor{gray!25} \textbf{1.5} & 5.0 & \textbf{3.0} & 5.0 & \cellcolor{gray!25} \textbf{1.5} & \textbf{7.0} \\
\bottomrule
\multicolumn{2}{l}{Total} & 5.143 & \cellcolor{gray!25} \textbf{2.714} & 5.357 & \textbf{3.643} & 5.143 & \textbf{3.214} & \textbf{2.786} \\
\bottomrule
\end{tabular}
\label{tbl:results-constr-multi}
\end{table*}